\documentclass[11pt,twoside]{article}
\usepackage{mathbbold}
\usepackage{amsmath, amssymb, mathrsfs, graphicx}

\def\scr{\mathscr}

\def\az{\alpha}  \def\bz{\beta}
    
    \def\fz{\varphi}
\def\gz{\gamma}  
\def\lz{\lambda}

\def\qz{\psi}
        
        \def\uz{\theta}
\def\vz{\varepsilon}

\def\ffz{\Phi}    \def\qqz{\Psi}
\def\ggz{\Gamma}

\def\qd{\quad}
\def\qqd{\qquad}

\newcommand{\mathsym}[1]{{}}

\def\scr{\mathscr}

\def\le{\leqslant}
\def\ge{\geqslant}

\font\cms=cmss9 scaled \magstep1

\newcommand{\nnb}{\nonumber}
\def\nnd{\noindent}

\def\thm{\nnd\bg{thm1}}
\def\crl{\nnd\bg{crl1}}
\def\lmm{\nnd\bg{lmm1}}
\def\prp{\nnd\bg{prp1}}

\def\xmp{\nnd\bg{xmp1}}

\def\hyp{\nnd\bg{hyp1}}
\def\dethm{\end{thm1}}
\def\decrl{\end{crl1}}
\def\delmm{\end{lmm1}}
\def\deprp{\end{prp1}}
\def\dexmp{\end{xmp1}}

\def\dehyp{\end{hyp1}}

\def\prf{\medskip \noindent {\bf Proof}. }
\def\deprf{\quad $\square$ \medskip}
\def\bg{\begin}
\def\be{\bg{equation}}
\def\de{\end{equation}}

\def\dear{\end{eqnarray}}
\def\lb{\label}
\def\ct{\cite}

\newcommand{\rf}[2]{[\ref{#1}; #2]}

\def\den{\end{enumerate}}

\def\ent{\text{\rm Ent${}_{\pi}$}}
\def\d{\text{\rm d}}

\begin{document}

\allowdisplaybreaks[4]
\thispagestyle{empty}
\renewcommand{\thefootnote}{\fnsymbol{footnote}}

\noindent {Acta Mathematica Sinica, English Series, 2012?}

\vspace*{.5in}
\begin{center}
{\bf\Large Bilateral Hardy-type Inequalities}
\vskip.15in {Mu-Fa Chen}
\end{center}
\begin{center} (Beijing Normal University, Beijing 100875, China)\\
\vskip.1in
May 2, 2012 \end{center}
\vskip.1in

\markboth{\sc Mu-Fa Chen}{\sc Bilateral Hardy-type inequalities}


\date{May 2, 2012}


\footnotetext{Received ??? 2012; accepted ??? 1012}
\footnotetext{{\it Key words and phases}.\quad
Hardy-type inequality, vanishing
at two endpoints, mean zero, splitting technique, normed linear space, Nash inequality, logarithmic Sobolev inequality.}

\bigskip

\begin{abstract}
This paper studies the Hardy-type inequalities on the intervals (may be infinite) with two weights, either vanishing
at two endpoints of the interval or having mean zero. For the first type of inequalities, in terms of new isoperimetric constants,
the factor of upper and lower bounds becomes smaller than the known ones.
The second type of the inequalities is motivated from probability theory and is new in the analytic context. The proofs are now rather elementary. Similar improvements are made for Nash inequality,
Sobolev-type inequality, and the logarithmic Sobolev inequality on the intervals.
\end{abstract}

\medskip

A large number of results on Hardy-type inequalities have been already collected and explored in the books \ct{kmp10} -- \ct{kp03}, \ct{mav},
and \ct{ok90}. This paper makes two additions. The first one is for the functions
vanishing at two endpoints of the interval. This type of inequalities was included in \ct{ok90}.
The contribution here is some improvement, not only on the isoperimetric constant but
also on the factor of the upper and lower bounds. The second addition is for the case where the functions have mean zero, which is motivated from a
probabilistic consideration and is not included in the books cited above. These two cases
are studied in the next two sections separately.
The main result in each case is stated as a theorem (Theorems \ref{t1-7} and \ref{t2-6}). Their extensions to more general
setup are presented as Theorems \ref{t1-9} and \ref{t2-9}.
As applications of the results or ideas developed in the first two section, in the third
section, we study the Nash inequality, the
Sobolev-type inequality, and the logarithmic Sobolev inequality.
The paper can be regarded as an extension
of the $L^2$-case studied in \ct{cmf10, cmf11}.

\section{The case of vanishing boundaries.}

Consider an interval $[-M, N]$ with $M, N\le \infty$. Certainly, here
$[-M, N]$ means $[-M, N)$ if $N=\infty$. This costs no confusion.
In this section, we study the case that the functions vanish
at two endpoints of the interval. That is the Hardy-type inequality:
\be \bigg(\int_{-M}^N |f|^q\d\mu\bigg)^{1/q}\!\!\le A
 \bigg(\int_{-M}^N |f'|^p\d\nu\bigg)^{1/p},\quad f(-M)\!=\!0\text{ and } f(N)\!=\!0,\lb{01}\de
where $f(N)=0$ for instance means that $\lim_{x\to\infty}f(x)=0$
if $N=\infty$.
Throughout the paper, all of the functions involved in the Hardy-type inequalities are assumed to be
absolutely continuous without mentioned time by time in what follows.
In the special case that $p=q=2$, the results in this and the next sections are
proved in \ct{cmf10, cmf11}) using much advanced methods. The present study is motivated from seeking for more direct proofs for the results.
At the moment, it is unclear how the capacitary technique used in
\ct{cmf10, cmf11} can be applied in the present general setup.
It may be helpful to the reader by studying the problem step by step to show how to find out the main result.
The study consists of five steps.
At each step, we have either a proposition
or a lemma. If one is in hurry, who may jump from here to the
main results, Theorems \ref{t1-9} and \ref{t2-9}.

Our first step is using the splitting technique (which we have used several times
before, cf. \ct{cmwf98}, \rf{cmf00}{Theorems 3.3 and 3.4}, and \ct{cmf05}). To do so,
fix $\uz\in (-M, N)$ and denote by $A_{\uz}^+$ and $A_{\uz}^-$, respectively, the optimal constant in the following inequalities.
$$\aligned
&\bigg(\int_{\uz}^N |f|^q\d\mu\bigg)^{1/q}\le A_{\uz}^+ \bigg(\int_{\uz}^N |f'|^p\d\nu\bigg)^{1/p},\qquad f(N)=0,\\
&\bigg(\int_{-M}^{\uz} |f|^q\d\mu\bigg)^{1/q}\le A_{\uz}^- \bigg(\int_{-M}^{\uz} |f'|^p\d\nu\bigg)^{1/p},\qquad f(-M)=0.
\endaligned$$
Clearly, these inequalities are different from (\ref{01}) since only one-side boundary
condition is endowed. Here and in what follows the superscript
 ``$-$'' means on the left-hand side of $\uz$
and ``$+$'' means on the right-hand side of $\uz$.

The next result shows that we can describe the optimal constant
$A$ in (\ref{01}) in terms of $A_{\uz}^{\pm}$ which are the optimal
constants on half-spaces with different boundary conditions.

\prp\lb{t1-1}{\cms For $1\le p\le q< \infty$, we have
$$2^{1/q-1/p} \sup_{\uz\in [-M, N]}\big(A_{\uz}^- \wedge A_{\uz}^+\big)\le A \le \inf_{\uz\in [-M, N]} \big(A_{\uz}^- \vee A_{\uz}^+\big),
$$
where $A_N^+=0$ and $A_{-M}^-=0$ by convention, $\az\wedge \bz=\min\{\az, \bz\}$,
and $\az\vee \bz=\max\{\az, \bz\}$.}
\deprp

\prf (a)
For each fixed $\uz\in [-M, N]$ and $f$ with $f(-M)=0$ and $f(N)=0$,
by the inequalities on the half-spaces, we have
$$\aligned
\int_{-M}^N |f'|^p\d\nu&=\int_{-M}^{\uz} |f'|^p\d\nu+\int_{\uz}^N |f'|^p\d\nu\\
&\ge  \big(A_{\uz}^-\big)^{-p}\bigg(\int_{-M}^{\uz} |f|^q\d\mu\bigg)^{p/q}+  \big(A_{\uz}^+\big)^{-p}\bigg(\int_{\uz}^N |f|^q\d\mu\bigg)^{p/q}\\
&\ge \Big[\big(A_{\uz}^-\big)^{-p}\wedge \big(A_{\uz}^+\big)^{-p}\Big]\bigg[\bigg(\int_{-M}^{\uz} |f|^q\d\mu\bigg)^{p/q}
  + \bigg(\int_{\uz}^N |f|^q\d\mu\bigg)^{p/q}\bigg]\\
&\ge \big(2^{(p/q-1)\vee 0}\big)^{-1}\Big[\big(A_{\uz}^-\big)^{-p}\wedge \big(A_{\uz}^+\big)^{-p}\Big] \bigg(\int_{-M}^{N} |f|^q\d\mu\bigg)^{p/q}\\
&\qquad\qquad\quad(\text{by $c_r$-inequality}).
\endaligned$$
Since $f$ is arbitrary, we have
$$A^p\le 2^{(p/q-1)\vee 0} \Big[\big(A_{\uz}^-\big)^{p}\vee \big(A_{\uz}^+\big)^{p}\Big].$$
Now, since $\uz$ is arbitrary, we obtain
$$A\le 2^{(1/q-1/p)\vee 0} \inf_{\uz\in [-M, N]}\big(A_{\uz}^- \vee A_{\uz}^+\big).$$
This conclusion holds for general $p, q\in [1, \infty)$.

(b) Again, fix $\uz$. Suppose for a moment that we can construct
two absolutely continuous functions $f_-$ and $f_+$ having the following properties:
$f_-(\!-\!M)$ $=0$, $f_-'(\uz)=0$, $f_-(\uz)>0$,
$$\int_{-M}^{\uz} |f_-|^q\d\mu=1,\qquad\text{and } \bigg(\int_{-M}^{\uz} |f_-'|^p\d\nu\bigg)^{1/p}< \big(A_{\uz}^-\big)^{-1}+\vz;$$
$f_+(N)=0$, $f_+'(\uz)=0$, $f_+(\uz)>0$,
$$\int_{\uz}^N |f_+|^q\d\mu=1,\qquad\text{and } \bigg(\int_{\uz}^N |f_+'|^p\d\nu\bigg)^{1/p}< \big(A_{\uz}^+\big)^{-1}+\vz.$$
Set $f=c f_-\mathbbold{1}_{[-M, \uz]} +f_+\mathbbold{1}_{(\uz, N]}$, where $c=f_+(\uz)/f_-(\uz)$. Then
$$\aligned
1+|c|^q=\int_{-M}^{\uz} &|cf_-|^q\d\mu+\int_{\uz}^N |f_+|^q\d\mu=\int_{-M}^{N} |f|^q\d\mu,\\
\int_{-M}^N |f'|^p\d\nu&=|c|^p \int_{-M}^{\uz} |f_-'|^p\d\nu+\int_{\uz}^N |f_+'|^p\d\nu\\
&\le |c|^p \Big(\big(A_{\uz}^-\big)^{-1}+\vz\Big)^p+ \Big(\big(A_{\uz}^+\big)^{-1}+\vz\Big)^p\\
&\le \Big(\big(A_{\uz}^-\big)^{-1}\vee \big(A_{\uz}^+\big)^{-1}+\vz\Big)^p (1+|c|^p).
\endaligned$$
Hence
$$\aligned
\bigg(\int_{-M}^N |f'|^p\d\nu\bigg)^{1/p}
&\le  \Big(\big(A_{\uz}^-\big)^{-1}\vee \big(A_{\uz}^+\big)^{-1}+\vz\Big) (1+|c|^p)^{1/p}\\
&\le 2^{1/p-1/q} \Big(\big(A_{\uz}^-\big)^{-1}\vee \big(A_{\uz}^+\big)^{-1}+\vz\Big) (1+|c|^q)^{1/q}\\
&\qquad\qquad\quad\text{(by Jensen's inequality requiring $q\ge p$)}\\
&=2^{1/p-1/q} \Big(\big(A_{\uz}^-\big)^{-1}\vee \big(A_{\uz}^+\big)^{-1}+\vz\Big) \bigg(\int_{-M}^{N} |f|^q\d\mu\bigg)^{1/q}.
\endaligned$$
Thus, whenever $q\ge p$, we have
$$A\ge 2^{1/q-1/p} \big(A_{\uz}^- \wedge A_{\uz}^+\big).$$
Therefore, we obtain
$$A\ge 2^{1/q-1/p} \sup_{\uz\in [-M, N]}\big(A_{\uz}^- \wedge A_{\uz}^+\big),
\qquad 1\le p\le q< \infty.$$
Combining this with (a), we arrive at the conclusion of the proposition.

(c) To complete the proof, it remains to construct the functions
$f_-$ and $f_+$ used in (b). For this, we need consider $f_-$ only by symmetry. The problem is the condition at $\uz$:
$f_-'(\uz)=0$ and $f_-(\uz)>0$.
The proof given below is modified from \rf{czz03}{Proof (ii) of Theorem 1.1}.
If necessary, by modifying $f_-$ properly on a sufficiently small neighborhood of $\uz$,
we can assume that $f_-'(\uz)=0$. The main point here
is to modify $f_-$ so that we also have $f_-(\uz)\ne 0$.
Otherwise, suppose that $f_-(\uz)= 0$.
Since $f_-$ is absolutely continuous, $f_-(-M)=0$ and $f_-(\uz)=0$, there exists $x_1\in (-M, \uz)$ such that
$|f_-(x_1)|=\sup_{x\in (-M, \uz)}|f_-(x)|$. Then $f_-(x_1)\ne 0$ (otherwise, $f\equiv 0$ which
contradicts with the norm 1 assumption). Let
$\tilde f_-=f_- \mathbbold{1}_{[-M, x_1)}+f_-(x_1)\mathbbold{1}_{[x_1, \uz]}$. Then
$\tilde f_-$ is absolutely continuous,
$$\aligned
&c^q:=\int_{-M}^{\uz}\big|\tilde f_-\big|^q \d \mu
\ge \int_{-M}^{\uz}\big|f_-\big|^q \d \mu=1,\\
&\bigg(\int_{-M}^{\uz}\big|\tilde f_-'\big|^p \d \nu\bigg)^{1/p}
\le \bigg(\int_{-M}^{\uz}\big|f_-'\big|^p \d \nu\bigg)^{1/p}
< \big(A_{\uz}^-\big)^{-1}+\vz.
\endaligned$$
Set $\bar f_-=c\tilde f_-$. Now it follows that
$$\bar f_-(-M)=0,\quad\bar f_-'(\uz)=0,\quad \bar f_-(\uz)\ne 0,
\quad \int_{-M}^{\uz}\big|\bar f_-\big|^q \d \mu=1,$$ and
$$\bigg(\int_{-M}^{\uz}\big|\bar f_-'\big|^p \d \nu\bigg)^{1/p}
=\frac{1}{c} \bigg(\int_{-M}^{\uz}\big|\tilde f_-'\big|^p \d \nu\bigg)^{1/p}
< \big(A_{\uz}^-\big)^{-1}+\vz.$$
Hence, we can replace $f_-$ by $\bar f_-$ when $f_-(\uz)=0$.
\deprf

Having Proposition \ref{t1-1} at hand, it is ready to write down some estimates of the optimal constant $A$ in (\ref{01}), as we did in \ct{cmf00, cmf05}, in terms of $B_{\uz}^{\pm}$ given below (cf. \rf{ok90}{Theorem 6.2} and \rf{mav}{\S 1.3, Theorem 3} in which the factor $ k_{q, p}$ may be different):
\begin{gather}
B_{\uz}^{\pm}\le A_{\uz}^{\pm}\le k_{q, p} B_{\uz}^{\pm},\qqd
1<p\le q<\infty,\lb{02}\\
B_{\uz}^+=\sup_{r\in (\uz,\, N)}\mu [\uz, r]^{1/q}
\bigg[\int_r^N \bigg(\frac{\d \nu^*}{\d x}\bigg)^{-1/(p-1)}\d x\bigg]^{(p-1)/p},\lb{03}\\
B_{\uz}^-=\sup_{r\in (-M,\, \uz)}\mu [r, \uz]^{1/q}\bigg[\int_{-M}^r \bigg(\frac{\d \nu^*}{\d x}\bigg)^{-1/(p-1)}\d x\bigg]^{(p-1)/p},\lb{04}
\end{gather}
where $\nu^*$ is the absolutely continuous part of $\nu$ and $k_{q, p}$
is a universal constant will be used often in this paper:
\be k_{q, p}=\bigg(1+\frac{q}{p'}\bigg)^{1/q}\bigg(1+\frac{p'}{q}\bigg)^{1/p'},\lb{04-1}\de where $p'$ is the conjugate number of $p$ and similarly for $q'$.
On the half-line when $q>p$, the constant is improved as follows:
$$
k_{q, p}\!=\!\Bigg[\frac{\Gamma\big(\frac{pq}{q-p}\big)}{\Gamma\big(\frac{q}{q-p}\big)\Gamma\big(\frac{p(q-1)}{q-p}\big)}\Bigg]^{1/p-1/q}
\!\!=\!\bigg[\frac{q-p}{pq B\big(\frac{q}{q-p},\, \frac{p(q-1)}{q-p}\big)}\Bigg]^{1/p-1/q},\qqd
q>p,
$$
where $\Gamma(x)$ and $B(x, y)=\Gamma(x)\Gamma(y)/\Gamma(x+y)$ are Gamma and Beta functions, respectively (cf. \rf{bg91}{Theorem 8}, \rf{mvm92}{Theorem 2}, and also \rf{kmp}{pages 45--47} for historical remarks). According to Lebesgue's decomposition theorem, each measure
$\nu$ can be decomposed into three parts:
$$\nu=\nu^*+\nu_{\text{\rm sing}}+\nu_{\text{\rm pp}},$$
where $\nu_{\text{\rm sing}}$ is the singular continuous part
and $\nu_{\text{\rm pp}}$ is the pure point part (a discrete measure).

We are now going to present some more explicit estimates.
To do so, we need the following simple result.

\lmm\lb{t1-2}{\cms For a given Borel measure $\mu$ and positive functions
$\fz$ and $\qz$ on $[-M, N]$, we have
$$\sup_{(x,\, y):\, x\le y}\frac{\mu[x, y]}{\fz(x)+\qz(y)}
\ge \sup_{\uz}\bigg\{\bigg[\sup_{x\le\uz}\frac{\mu[x, \uz]}{\fz(x)}\bigg]\bigwedge \bigg[\sup_{y\ge\uz}\frac{\mu(\uz, y]}{\qz(y)}\bigg]\bigg\}.$$}
\delmm

\prf For fixed $x\le y$ and $(x, y)\ni \uz$, we have by proportional property that
$$\frac{\mu[x,\, y]}{\fz(x)+\qz(y)}
= \frac{\mu[x, \uz]+\mu(\uz, y]}{\fz(x)+\qz(y)}
\ge \frac{\mu[x, \uz]}{\fz(x)}\bigwedge \frac{\mu(\uz, y]}{\qz(y)}$$
and furthermore
$$\frac{\mu[x, y]}{\fz(x)+\qz(y)}
\ge\sup_{\uz\in [x,\, y]}\bigg\{ \frac{\mu[x, \uz]}{\fz(x)}\bigwedge \frac{\mu(\uz, y]}{\qz(y)}\bigg\}.$$
Thus,
$$
\aligned
\sup_{x\le y}\frac{\mu[x, y]}{\fz(x)+\qz(y)}
&\ge\sup_{x\le y}\sup_{\uz\in [x, y]}\bigg\{ \frac{\mu[x, \uz]}{\fz(x)}\bigwedge \frac{\mu(\uz, y]}{\qz(y)}\bigg\}\\
&=\sup_{\uz}\,\sup_{[x,\, y]\ni\uz}\big\{\cdot\cdot\cdot\big\}\\
&=\sup_{\uz}\bigg\{\bigg[\sup_{x\le\uz}\frac{\mu[x, \uz]}{\fz(x)}\bigg]\bigwedge \bigg[\sup_{y\ge\uz}\frac{\mu(\uz, y]}{\qz(y)}\bigg]\bigg\}\endaligned$$
as required.
\deprf

It is remarkable that we do not have an expected dual result of
the above lemma. At beginning, we do have the dual
$$\frac{\mu[x, y]}{\fz(x)+\qz(y)}
\le \frac{\mu[x, \uz]+\mu(\uz, y]}{\fz(x)+\qz(y)}
\le \frac{\mu[x, \uz]}{\fz(x)}\bigvee \frac{\mu(\uz, y]}{\qz(y)}.$$
Hence
$$\sup_{x\le\uz\le y}
\frac{\mu[x, y]}{\fz(x)+\qz(y)}
\le \bigg[\sup_{x\le\uz}\frac{\mu[x, \uz]}{\fz(x)}\bigg]\bigvee \bigg[\sup_{y\ge\uz}\frac{\mu(\uz, y]}{\qz(y)}\bigg]$$
and furthermore
$$\inf_{\uz}\sup_{x\le\uz\le y}
\frac{\mu[x, y]}{\fz(x)+\qz(y)}
\le \inf_{\uz}\bigg[\sup_{x\le\uz}\frac{\mu[x, \uz]}{\fz(x)}\bigg]\bigvee \bigg[\sup_{y\ge\uz}\frac{\mu(\uz, y]}{\qz(y)}\bigg].$$
Clearly, this is somehow a dual of Lemma \ref{t1-2} but it is still a distance to what we expect:
$$\sup_{x\le y}\frac{\mu[x, y]}{\fz(x)+\qz(y)}
\le \inf_{\uz}\bigg[\sup_{x\le\uz}\frac{\mu[x, \uz]}{\fz(x)}\bigg]\bigvee \bigg[\sup_{y\ge\uz}\frac{\mu(\uz, y]}{\qz(y)}\bigg].$$
Alternatively, let $\bar\uz$ satisfy
$$\sup_{x\le\bar\uz}\frac{\mu[x, \bar\uz]}{\fz(x)}= \sup_{y\ge\bar\uz}\frac{\mu(\bar\uz, y]}{\qz(y)}.$$
Then we have
\be\sup_{x\le{\bar\uz}\le y}\frac{\mu[x, y]}{\fz(x)+\qz(y)}
\le \sup_{x\le\bar\uz}\frac{\mu[x, \bar\uz]}{\fz(x)}=
\bigg[\sup_{x\le\bar\uz}\frac{\mu[x, \bar\uz]}{\fz(x)}\bigg]\bigvee \bigg[\sup_{y\ge\bar\uz}\frac{\mu(\bar\uz, y]}{\qz(y)}\bigg].\lb{05-1}\de
Very often, the right-hand side coincides with
$$\inf_{\uz}\bigg\{\bigg[\sup_{x\le\uz}\frac{\mu[x, \uz]}{\fz(x)}\bigg]\bigvee \bigg[\sup_{y\ge\uz}\frac{\mu(\uz, y]}{\qz(y)}\bigg]\bigg\},$$
but one can not remove $\bar\uz$ from the left-hand side and
keep the inequality.

Throughout this paper, we mainly restrict ourselves to the case that $1<p\le q<\infty$. The limit case that
either $p=1$ or $q=\infty$ are easier and so are omitted here. For simplicity,
throughout this paper, we set
$$h(x)= \bigg(\frac{\d \nu^*}{\d x}\bigg)^{-1/(p-1)},\qqd {\hat\nu}(\d x)=h(x)\d x.$$
Clearly, $h$ and $\hat\nu$ depend on $p>1$. The measure $\hat\nu$ comes, but different, from $\nu$.
In what follows, almost every estimate is expressed by using the pair $(\mu, \hat\nu)$ but not
$(\mu, \nu)$. Besides, we may assume that
\be {\hat\nu}(-M, N):=\int_{-M}^N  h=\int_{-M}^N \bigg(\frac{\d \nu^*}{\d x}\bigg)^{-1/(p-1)}\d x<\infty.\lb{05}\de
This technical assumption can often be avoided by replacing
$\d \nu^*/\d x$ with $\d \nu^*/\d x+\vz \exp\big[(p-1) x^{2}\big]$ and then passing to the limit as $\vz\downarrow 0$. Alternatively, one may start at $M, N<\infty$, replace $\d \nu^*/\d x$ with $\d \nu^*/\d x+\vz$. Then pass to the limit as $\vz\downarrow 0$, and then as $M, N\to\infty$ if necessary. In parallel, without loss of generality, we can also assume that $\mu$ is positive on each subinterval.

Next, define a constant $B^*$ by
\be \!\big({B^*}^q\big)^{-1}\!=\!\!
\inf_{-M\le x\le y\le N}\Big[{\hat\nu}[-M, x]^{-\frac{q(p-1)}{p}}+{\hat\nu}[y, N]^{-\frac{q(p-1)}{p}}\Big]
\mu[ x,  y]^{-1}.\lb{06}\de
Let us now discuss the boundary condition in the definition of ${B^*}$ above (or $B_*$ below),  when $M=\infty$, here $x=-M$ means that $x\to -\infty$:
$$\aligned
\varlimsup_{x\to -\infty}&\mu[x, y]^{1/q}\Big[{\hat\nu}[-M, x]^{-\frac{q(p-1)}{p}}+{\hat\nu}[y, N]^{-\frac{q(p-1)}{p}}\Big]^{-1/q}\\
&=\varlimsup_{x\to -\infty}\mu[x, y]^{1/q}\,{\hat\nu}[-M, x]^{(p-1)/p}\endaligned$$
which is the type $\infty\cdot 0$ of limit provided $\mu[-\infty, y]=\infty$.
Otherwise, the limit is zero and so the boundary $-M$ can be ignored in computing $B^*$.
When $M=\infty=N$, we need to compute the iterated limit only.
To which, the main reason is that the optimal constant $A$ is
increasing as either $N\uparrow$ or $-M\downarrow$. Hence, the
general case can be regarded as the limit of finite $M$ and $N$.
In other words, we do not need to consider the other types of
double limits as $M, N\to\infty$.

Here is our upper estimate.

\lmm\lb{t1-3} {\cms Let $\mu_{\text{\rm pp}}=0$. Then for $1<p\le q<\infty$, we have
$A\le k_{q, p} B^*,$
where $B^*$ is defined by (\ref{06})}.
\delmm

\prf As mentioned above, without loss of generality, we can assume (\ref{05}). Rewrite ${B^*}^q$ as
$${B^*}^q=\sup_{x\le y}
\frac{\mu[x,  \uz]+ \mu(\uz, y]}{
{\hat\nu}[-M, x]^{-q(p-1)/p}+{\hat\nu}[y, N]^{-q(p-1)/p}}.$$
As an application of Lemma \ref{t1-2}, we have
$$\aligned
{B^*}^q&\ge
\sup_{\uz}\bigg\{\bigg[\sup_{x\le\uz}\frac{\mu[x, \uz]}{{\hat\nu}[-M, x]^{-q(p-1)/p}}\bigg]\bigwedge \bigg[\sup_{y\ge\uz}\frac{\mu(\uz, y]}{{\hat\nu}[y, N]^{-q(p-1)/p}}\bigg]\bigg\}\\
&=\sup_{\uz}\big[B_{\uz}^-\wedge B_{\uz}^+\big]^q.
\endaligned$$
Here in the last step, we have used the condition $\mu_{\text{\rm pp}}=0$. Since we can represent
$$\mu_{\text{\rm sing}}[x, \uz]=\mu_{\text{\rm sing}}[-M, \uz]- \mu_{\text{\rm sing}}[-M, x) $$
(the last term is continuous in $x$) and
$\mu_{\text{\rm pp}}=0$, it follows that the function $\mu[x, \uz]$ is continuous in $x$ and $\uz$.
By choosing $\bar\uz$ such that $B_{\bar\uz}^-= B_{\bar\uz}^+$, it follows that $B_{\bar\uz}^-\le B^*$ (just proved) and furthermore
$$\aligned
A&\le \inf_{\uz\in [-M, N]}\big(A_{\uz}^-\vee A_{\uz}^+\big)\qd(\text{by Proposition \ref{t1-1}})\\
&\le k_{q, p} \inf_{\uz\in [-M, N]}\big(B_{\uz}^-\vee B_{\uz}^+\big)\qd(\text{by (\ref{02})})\\
&\le k_{q, p} B_{\bar\uz}^-\qd(\text{by definition of $\bar\uz$})\\
&\le k_{q, p} {B^*}.\qqd\square
\endaligned$$
\smallskip

Note that the parameter $\uz$ is used temporary in the proof above. Thus, the splitting procedure
is a bridge to go to the upper estimate but our final result does not depend on the splitting points. The simple technique used in
proving the upper estimate above is in common, and will be used several times later (Lemma \ref{t2-3} and Theorem \ref{t3-2}).

Before moving further, let us discuss the technical assumption that $\mu_{\text{\rm pp}}=0$.  In the study of $B_{\uz}^{\pm}$ for half-spaces, one first handles with the case that
$\mu\ll\d x$ and $\nu\ll \d x$ and then removes this restriction by the following technique. Without loss of
generality, assume that $M, N<\infty$. Besides, we may also assume that
$f'\ge 0$ in the study of the upper estimate. The idea is to use an approximating procedure (cf. \ct{mb72} or \rf{mav}{page 45}). Note that
$$\bigg(\int_{-M}^N f^q\d \mu\bigg)^{1/q}
=\bigg(\int_{-M}^N \mu [x, N]\, \d f (x)^q\bigg)^{1/q}.$$
Now, we can approximate $\mu[x, N]$ by a sequence of absolutely
continuous, decreasing functions $\{g_n\}$ having the property: $g_n\le \mu[\cdot, N]$ for every $n$; as $n\to\infty$, $g_n(x)$ converges to $\mu[x, N]$ for almost all $x$.
Thus, we can first replace $\mu[x, N]$ by absolutely continuous
$g_n$ and then pass to the limit as $n\to\infty$. Actually, now
 $\mu[\cdot, N]$ consists of three parts: the absolutely continuous part, the singular continuous one plus a step function. Each of them is decreasing. There is nothing to do about the absolutely continuous part. The singular decreasing continuous function can be approximated from below by decreasing step functions. Furthermore, each of the step functions can be approximated from below almost everywhere by absolutely continuous decreasing functions.  The new difficulty arises: even though we have the control $g_n(x)\le \mu[x, N]$ for all $x$, but we still do not know how to construct
a sequence $\{g_n\}$ as above having the control $(0\le)\,g_n(x)-g_n(y)\le \mu[x, y]$ for every pair $\{x, y\}$ with $x<y$ and each $n$.

For the lower bound of $A$, the dual proof of Lemma \ref{t1-2} does not work well as remarked below Lemma \ref{t1-2}. More precisely, what we can obtain by Proposition \ref{t1-1} and (\ref{05-1}) is as follows.
\be A\ge 2^{1/q-1/p}\sup_{x\le{\bar\uz}\le y}\bigg\{
\mu[x, y]^{1/q}\Big[{\hat\nu}[-M, x]^{-\frac{q(p-1)}{p}}+{\hat\nu}[y, N]^{-\frac{q(p-1)}{p}}\Big]^{-1/q}\bigg\},\de
where $\bar\uz$ is the solution of the equation $B_{\uz}^-=B_{\uz}^+$. The result is less satisfactory since $\bar\uz$ (unknown explicitly) is included. Fortunately, there
is a direct technique (cf. \rf{cmf10}{Proof (b) of Theorem 8.2})
to handle with the lower estimate.

\lmm\lb{t1-4} {\cms For $1<p,\, q<\infty$, we have
$$A\ge \sup_{-M\le x\le y\le N}\bigg\{\mu[x, y]^{1/q}\Big({\hat\nu}[-M, x]^{1-p}+{\hat\nu}[y, N]^{1-p}\Big)^{-1/p}\bigg\}=:B_*.$$
}
\delmm

\prf
Given $m', m, \uz, n, n'\in [-M, N]$ with $m'<m<\uz<n<n'$, define
$$
f(x)=\gz\mathbbold{1}_{\{m'\le x\le \uz\}} {\hat\nu}[m', x\wedge m]+
\mathbbold{1}_{\{\uz< x\le n'\}} {\hat\nu}[x\vee n, n']$$
where
$\gz= {\hat\nu}[n, n']/ {\hat\nu}[m', m].$
Clearly, $f$ is absolutely continuous. We have
$$\bigg(\int_{-M}^N |f|^q\d\mu\bigg)^{1/q}\!\!
=\!\bigg(\int_{m'}^{n'} |f|^q\d\mu\bigg)^{1/q}\!\!
\ge\! \bigg(\int_{m}^{n} |f|^q\d\mu\bigg)^{1/q}
\!\!=\mu[m, n]^{1/q} \,{\hat\nu}[n, n']$$
and
$$\aligned
\bigg(\int_{-M}^N |f'|^p\d\nu\bigg)^{1/p}
&=\bigg(\gz^p \int_{m'}^m h^p\d\nu+\int_n^{n'} h^p\d\nu\bigg)^{1/p}\\
&=\big(\gz^p \,{\hat\nu}[m', m]+{\hat\nu}[n, n']\big)^{1/p}.
\endaligned$$
Here in the last step, we have ignored the singular part of $\nu$
since the original inequality is equivalent to the one having
$\nu=\nu^*$. To see this, simply set $f'=0$ on the singular part of $\nu$. Thus, the optimal constant $A$ satisfies
$$A\ge \big({\mu[m, n]^{1/q}\, {\hat\nu}[n, n']}\big){\big(\gz^p \, {\hat\nu}[m', m]+{\hat\nu}[n, n']\big)^{-1/p}}.$$
But
$$
\aligned
&\big(\gz^p\, {\hat\nu}[m', m]+{\hat\nu}[n, n']\big){\hat\nu}[n, n']^{-p} \\
&=\big\{{\hat\nu}[n, n']^p\, {\hat\nu}[m', m]^{-p}\,{\hat\nu}[m', m]
+{\hat\nu}[n, n']\big\}\,{\hat\nu}[n, n']^{-p}\\
&={\hat\nu}[m', m]^{1-p}
+{\hat\nu}[n, n']^{1-p},
\endaligned$$
it follows that
$$A\ge \mu[m, n]^{1/q}\big({\hat\nu}[m', m]^{1-p}
+{\hat\nu}[n, n']^{1-p}\big)^{-1/p}.$$
Let $m'\downarrow -M$, $n'\uparrow N$ and then make supremum with respect to
$m=x\le y=n$. We get the required assertion.
\deprf

On the comparison of $B_*$ and $B^*$, it is obvious that $B^*=B_*$ if $p=q$. In general, we have the following result.

\lmm\lb{t1-4-1}{\cms Let $q\ge p$. Then we have $B_*\le B^*\le 2^{1/p-1/q}B_*$.}
\delmm

\prf Simply apply the $c_r$-inequality:
$$(\az+\bz)^r\le 2^{(r-1)\vee 0}(\az^r+\bz^r).$$

(a) Set
$$\az={\hat\nu}[-M, x]^{\frac{q(1-p)}{p}},\qd
\bz={\hat\nu}[y, N]^{\frac{q(1-p)}{p}},\qd
r=\frac p q \in (0, 1].$$
It follows that
$$(\az+\bz)^{1/q}\le \big(\az^{p/q}+\bz^{p/q}\big)^{1/p},$$
and then $B^*\ge B_*$.

(b) Set
$$\az={\hat\nu}[-M, x]^{1-p},\qd
\bz={\hat\nu}[y, N]^{1-p},\qd
r=\frac q p \ge 1.$$
We have
$$(\az+\bz)^{1/p}\le 2^{1/p-1/q}\big(\az^{q/p}+\bz^{q/p}\big)^{1/p},$$
and then $B_*\ge 2^{1/q-1/p}B^*$. Certainly, in this case
the assertion can also be deduced by Jensen's inequality.
\deprf

We mention that even though their supremums are equivalent but
in the proofs of Lemmas \ref{t1-3} and \ref{t1-4}, the expressions of $B^*$ and $B_*$ are not exchangeable, because $B_*$ does not own the homogeneous of that of $B^*$.

We are now ready to state
our first main result.

\thm\lb{t1-7}{\cms The optimal constant $A$ in the Hardy-type inequality (\ref{01}) satisfies

(1)\; $A \le k_{q, p} B^*$ for $1< p\le q< \infty$ once
$\mu_{\text{\rm pp}}=0$, where $k_{q, p}$ is defined by (\ref{04-1}), and

(2)\; $A\ge B_*$ for $1<p, q<\infty$, where
$$\aligned
B^*&=
\sup_{-M\le x\le y\le N}\bigg\{\mu[x, y]^{1/q}\Big({\hat\nu}[-M, x]^{q(1-p)/p}+{\hat\nu}[y, N]^{q(1-p)/p}\Big)^{-1/q}\bigg\} ,\\
B_*&=\sup_{-M\le x\le y\le N}\bigg\{\mu[x, y]^{1/q}\Big({\hat\nu}[-M, x]^{1-p}+{\hat\nu}[y, N]^{1-p}\Big)^{-1/p}\bigg\}.
\endaligned$$
Moreover, we have $B_*\le B^*\le 2^{1/p-1/q} B_*$ when $q\ge p$.}
\dethm

\prf The conclusions are combination of Lemmas \ref{t1-3}--\ref{t1-4-1}.
\deprf

It is interesting to have a look at the factor $k_{q, p}$ in Theorem \ref{t1-7}. When $p=q$, the factor becomes
$$\bigg(\frac{q}{q-1}\bigg)^{(q-1)/q}q^{1/q}.$$
On $(1, \infty)$, it is unimodal having maximum 2 at $q=2$ and
decreases to 1 as $q\to 1$ or $\infty$. More generally, the rough
ratio of the upper and lower bounds is no more than
$$\bigg(1+\frac{q}{p'}\bigg)^{1/q}\bigg(1+\frac{p'}{q}\bigg)^{1/p'}2^{1/p-1/q}$$
which is again $\le 2$ (for every $q\ge p$), having equality sign iff $p=q=2$.

The study on the inequality (\ref{01}) was began by P. Gurka in an unpublished paper using a common constant
$$\aligned
\widetilde B&=\sup_{-M\le x\le y\le N}\bigg\{\mu[x, y]^{\frac 1 q}\Big({\hat\nu}[-M, x]^{1-p}
\bigvee {\hat\nu}[y, N]^{1-p}\Big)^{-\frac 1 p}\bigg\}\\
&=\sup_{-M\le x\le y\le N}\Big\{\mu[x, y]^{\frac 1 q}\Big({\hat\nu}[-M, x]^{\frac{p-1}{p}}
\bigwedge {\hat\nu}[y, N]^{\frac{p-1}{p}}\Big)\Big\}
\endaligned$$
and having a universal factor 8. This $\widetilde B$ is closely related to $B_*$: replacing ``$+$'' with ``$\vee$'', we obtain $\widetilde B$ from $B_*$. Gurka's result was then improved in \rf{ok90}{Theorem 8.2} with a smaller factor
(unexplicit one $\approx 4.71$ and explicit one $=2\sqrt{6}$ in the case of $p=q=2$). Note that using the inequalities
$$\az\vee\bz\le \az+\bz\le 2(\az\vee \bz),$$
from Theorem \ref{t1-7}, it follows that we have lower and upper
bounds replacing $B^*$ and $B_*$ by the same $\widetilde B$ with
an additional factor $2^{-1/p}$ for the lower estimate.
Then the factor becomes $2\sqrt{2}$ in the case of $p=q=2$.
Replacing $\az\vee\bz$ with $\az+\bz$ is an essential difference of the present paper from the previous ones in the bilateral
situation. Besides, the inequality (\ref{01}) was also
proved in \rf{ok90}{Theorem 8.8} with a common constant
$$\widetilde B\!=\!\!\inf_{-M\le \uz\le N}\!
\bigg\{
\Big[\!\sup_{-M\le x\le \uz}\mu[x, \uz]^{\frac 1 q}\,
{\hat\nu}[-M, x]^{\frac{p-1}{p}}\Big]
\!\bigvee\! \Big[\sup_{\uz\le y\le N}\mu[\uz, y]^{\frac 1 q}\,{\hat\nu}[y, N]^{\frac{p-1}{p}}\Big]\!\bigg\}$$
having a factor
$$2^{\frac 1 p} \bigg(1+ \frac{q}{p'}\bigg)^{\frac 1 q}
\bigg(1+\frac{p'}{q}\bigg)^{\frac{1}{p'}}$$
which has an additional factor is $2^{1/p}$ than (\ref{04-1}).
The last result is related to our splitting technique. All of
these results use the assumption that $\mu\ll \d x$ and $\nu\ll \d x$.

Before moving further, we want to describe $B_*$ and $B_*$ more carefully. It also leads some quantities which are easier in
practical computations. For this, we need some preparation.
Assume that (\ref{05}) holds.
For each $x\in (-M, N)$, let $y(x)$ be the unique solution of the equation
$${\hat\nu}[-M, x]= {\hat\nu}[y, N].$$
Next, let $m(\hat\nu)$ be a solution to the equation
$$y(x)=x, \qqd x\in (-M, N).$$ Thus,  $m(\hat\nu)$ is actually the median of the measure $\hat\nu$ (but not $\nu$):
$${\hat\nu}[-M, m] = {\hat\nu}[m, N].$$
Set
$$\aligned
H_{\mu, \nu}(x, y)=&\mu[x, y]^{1/q}\Big[{\hat\nu}[-M, x]^{1-p}
+{\hat\nu}[y, N]^{1-p}\Big]^{-1/p},\\
&\qqd\qqd\qqd\qqd\qqd -M\le x\le y\le N.\endaligned$$
Define
\be \! H^o=2^{-1/p}\!\!\sup_{x\in (-M,\, m(\hat\nu)]}\mu[x, y(x)]^{1/q}\,
{\hat\nu}[-M, x]^{(p-1)/p}\!.\lb{07}\de
Denote by $\ggz$ be the limiting points of $H_{\mu, \nu}(x, y)$ as $\mu[y, N]=\infty$
or $\mu[-M, x]=\infty$, as well as the iterated limits if $\mu[-M, N]=\infty$
when $M=\infty=N$. Set
\be H^{\partial}=
\begin{cases}
\sup \{\gz: \gz\in \ggz\}\qd &\text{if } \ggz\ne\emptyset\\
0 \qd &\text{if } \ggz=\emptyset.
\end{cases}\lb{08}\de
Clearly, $H^{\partial}=0$ if $M, N<\infty$.

We are now ready to describe $B^*$ and $B_*$ in terms of $H^o$ and $H^{\partial}$.

\lmm\lb{t1-5} {\cms Let (\ref{05}) hold. Then we have
$$H^o\vee H^{\partial}\le B_* \le \big(2^{1/p} H^o\big)\vee H^{\partial}.
$$}\delmm

\prf Rewrite $H$ as
$$H_{\mu, \nu}(x, y)=\left[\frac{{\hat\nu}[-M, x]^{1-p}
+{\hat\nu}[y, N]^{1-p}}{\mu[x, y]^{p/q}}\right]^{-1/p}.$$
under (\ref{05}), because for finite $x$ and $y$ with $x\le y$, we have
$$\frac{{\hat\nu}[-M, x]^{1-p}
+{\hat\nu}[y, N]^{1-p}}{\mu[x, y]^{p/q}}
\ge \frac{2}{{\mu[x, y]^{p/q}}}\bigg[{\hat\nu}[-M, x]\,
{\hat\nu}[y, N]\bigg]^{(1-p)/2},$$
and the equality sign holds iff ${\hat\nu}[-M, x]= {\hat\nu}[y, N]$
which gives us the solution $y(x)$. Thus, we obtain
\be \inf_{x\le y}\frac{{\hat\nu}[-M, x]^{1-p}
+{\hat\nu}[y, N]^{1-p}}{\mu[x, y]^{p/q}}
\le 2\inf_{x\le m(\hat\nu)} \frac{{\hat\nu}[-M, x]^{1-p}
}{\mu[x, y(x)]^{p/q}},\lb{09}\de
since $\{(x, y(x)): x\le m(\hat\nu)\}\subset \{(x, y): x, y\in (-M, N)\}$.
This gives us a lower bound of the supremum of $H_{\mu, \nu}(x, y)$
over the set $\{(x, y): x<y,\, \mu[x, y]<\infty\}$. Next, we have
$$\aligned
\inf_{x\le y}\frac{{\hat\nu}[-M, x]^{1-p}
+{\hat\nu}[y, N]^{1-p}}{\mu[x, y]^{p/q}}
&=\inf_{x\le y}\bigg\{\frac{{\hat\nu}[-M, x]^{1-p}
}{\mu[x, y]^{p/q}} +
\frac{{\hat\nu}[y, N]^{1-p}}{\mu[x, y]^{p/q}}
\bigg\}\\
&\ge \inf_{x\le y}\bigg\{\frac{{\hat\nu}[-M, x]^{1-p}
}{\mu[x, y]^{p/q}}
\bigvee
\frac{{\hat\nu}[y, N]^{1-p}}{\mu[x, y]^{p/q}}
\bigg\}\\
&=:\xi.\endaligned$$
Without loss of generality, assume that $M, N<\infty$.
Because of the continuity of the involved functions,
the minimum $\xi$ can be achieved at some pair $(x_0, y_0)$.
We now prove that $(x_0, y_0)$ should be located at the surface
where the two terms in the last $\{\cdots\}$ are equal. Otherwise,
without loss of generality, assume that
$$\vz:={\mu[x_0, y_0]^{-p/q}}\,{{\hat\nu}[-M, x_0]^{1-p}
}
-{\mu[x_0, y_0]^{-p/q}}\,{{\hat\nu}[y_0, N]^{1-p}}>0.$$
Let $\bar y>y_0$ be sufficiently close to $y_0$. Then we have
$${\mu[x_0, \bar y]^{-p/q}}\,{{\hat\nu}[-M, x_0]^{1-p}
}<{\mu[x_0, y_0]^{-p/q}}\,{{\hat\nu}[-M, x_0]^{1-p}
}$$
(here we have used the preassumption that $\mu$ is positive on each subinterval) and
$${\mu[x_0, \bar y]^{-p/q}}\,{{\hat\nu}[\bar y, N]^{1-p}
}
<{\mu[x_0, y_0]^{-p/q}}\,{{\hat\nu}[y_0, N]^{1-p}
}+ {\vz}/{2},$$
due to the continuity of the involved functions.
We have thus obtained a pair $(x_0, \bar y)$ with $x_0< \bar y$ such that
$$\Big\{\Big[{\mu[x_0, \bar y]^{-p/q}}\,{{\hat\nu}[-M, x_0]^{1-p}
}\Big]
\bigvee
\Big[{\mu[x_0, \bar y]^{-p/q}}\,{{\hat\nu}[\bar y, N]^{1-p}}\Big]
\Big\}<\xi.$$
This is a contradiction to the minimum property of $\xi$. Therefore, we obtain
\be \inf_{x\le y}{\mu[x, y]^{-p/q}}\Big[{{\hat\nu}[-M, x]^{1-p}
+{\hat\nu}[y, N]^{1-p}}\Big]
\ge \inf_{x\le m(\hat\nu)} \frac{{\hat\nu}[-M, x]^{1-p}
}{\mu[x, y(x)]^{p/q}}.\lb{10}\de
From this, we obtain a upper bound of the supremum of $H_{\mu, \nu}(x, y)$
over the set $\{(x, y): x<y,\, \mu[x, y]<\infty\}$ in terms of $H^o$ up to a factor $2^{-1/p}$.
In other words, we  have worked out the case that the supremum is achieved inside of the interval.
In general, it may be achieved at the $\infty$-boundaries (at which $\mu[y, N]=\infty$ or $\mu[-M, x]=\infty$).
This leads to the boundary condition $H^{\partial}$, when one of $M$ and $N$ is infinite.
Combining these two parts together, we get the estimates of  $B_*$ under (\ref{05}).
\deprf

An easier way to understand what was going on in the last proof
is look at the following simple example. Consider functions
$f(x)=2x$ and $g(x)=3-x$ on $[0, 2]$. They intersects uniquely at the point $x^*=1$. Then we have
$$2\sqrt{f(x^*)g(x^*)}=4
>\! \inf_{x\in [0, 2]}\big(f(x)+g(x)\big)=3
>\! \inf_{x\in [0, 2]}\big[f(x)\vee g(x)\big]\!=\!
f(x^*)=2.$$
If we rewrite the first term as $2\inf_{x\in [0, 2]}\big(f(x)\vee g(x)\big)$, then it becomes obvious that
the middle term can be bounded by the first and the last ones. Certainly, the bounds are usually not sharp.

\lmm\lb{t1-6}{\cms Let (\ref{05}) hold. Then we have
$$\big(2^{1/p-1/q}H^o\big)\vee H^{\partial}\le B^*\le \big(2^{1/p}H^o\big)\vee H^{\partial}.$$}
\delmm

\prf Note that the difference of $B^*$ and $B_*$ is only the
summation terms. If one of the terms in the sum is ignored, then
the remaining terms coincide with each other. Thus, the boundary condition $H^{\partial}$ is the same for $B^*$ and $B_*$.
When $M, N<\infty$, the proof of the comparison of $B^*$ with $H^o$ is similar to the last one.
\deprf

We remark that in the degenerated case that (\ref{05}) does not hold, say ${\hat\nu}[y, N]$ $=\infty$, then we have obviously that $B_*=B^*$.

As a combination of the last two lemmas, we obtain the following
simple criterion.

\crl{\cms The Hardy-type inequality (\ref{01}) holds iff $H^o\vee H^{\partial}<\infty$.}
\decrl

\prf When (\ref{05}) holds, the assertion follows from the last two lemmas. Note that if $\hat\nu [y, N]=\infty$ for instance, we have
$B^*=B_*$ and so the assertion is described by $H^{\partial}$
only. \deprf

We now extend Theorem \ref{t1-7} to a more general setup which is mainly used
in interpolation of $L^p$-spaces. For this, we need a class of normed linear spaces
$(\mathbb B, \|\cdot\|_{\mathbb B},\, \mu)$ consisting of real Borel
measurable functions on a measurable space $(X, {\scr X},\, \mu)$.
We now modify the hypotheses on the normed linear spaces given in
\rf{cmf05}{Chapter 7} as follows.
\hyp\lb{t1-8}{\cms
\begin{itemize}\setlength{\itemsep}{-0.8ex}
\item [(H1)] In the case that $\mu(X)=\infty$, $\mathbbold{1}_K\in {\mathbb B}$ for all
   compact $K$. Otherwise,  $1\in {\mathbb B}$.
\item [(H2)] If $h\in \mathbb B$ and  $|f|\le h$, then  $f\in \mathbb B$.
\item [(H3)] $\|f\|_{\mathbb B}=\sup_{g\in {\scr G}} \int_X |f| g
\d\mu$,
\end{itemize}
}\dehyp
where ${\scr G}$, to be specified case by case, is a class of
nonnegative $\scr X $-measurable functions. A typical example
is ${\scr G}=\{\mathbbold{1}\}$ and then ${\mathbb
B}=L^1(\mu)$. In what follows, the measure space
$(X, {\scr X},\, \mu)$ is fixed to be $\big([-M, N], {\scr B}([-M, N]), \mu\big)$. We often use the dual representation of the norm
$$\bigg(\int_{-M}^N |f|^r\d\mu \bigg)^{1/r}=\sup_{g\in\, \text{The unit ball in $L^{r'}(\mu)$}}\int_{-M}^N |f|g\d\mu,$$
where $r'$ is the conjugate number of $r\,(\ge 1)$.
Throughout this paper, we assume (H1)--(H3) for
$(\mathbb B, \|\cdot\|_{\mathbb B},\, \mu)$ without mentioned again.

For simplicity, we write the $L^p$-norm with respect to $\mu$
as $\|\cdot\|_{\mu, p}$. If necessary, we also write $\|\cdot\|_{\az, \bz; \mu, p}$
to indicate the interval $[\az, \bz]$.

\thm\lb{t1-9}{\cms Let ${\scr G}$ satisfy Hypotheses \ref{t1-8} and consider the Hardy-type inequality
$$\|f^q\|_{\mathbb B}^{1/q}\le A_{\mathbb B}\|f'\|_{\nu, p},\qqd f(-M)=0\text{\cms\;and } f(N)=0.$$
\begin{itemize}\setlength{\itemsep}{-0.8ex}
\item[(1)] Then the optimal constant $A_{\mathbb B}$ satisfies
$$A_{\mathbb B}\le k_{q, p}\, B_{\mathbb B}^*\qd \text{\cms\; for $1<p\le q<\infty$ once  $\mu_{\text{\rm pp}}=0$, and}$$
\item[(2)] $A_{\mathbb B}\ge B_{\mathbb B*}$ for $1<p, q<\infty$, where
$$\aligned
&B_{\mathbb B}^*=\sup_{-M\le x\le y\le N} \bigg\{\big\|\mathbbold{1}_{ [x, y]}\big\|_{\mathbb B}^{\frac 1 q}\Big({\hat\nu}[-M, x]^{\frac{q(1-p)}{p}}+{\hat\nu}[y, N]^{\frac{q(1-p)}{p}}\Big)^{-\frac 1 q}\bigg\},\\
&B_{\mathbb B*}=\sup_{-M\le x\le y\le N} \bigg\{\big\|\mathbbold{1}_{ [x, y]}\big\|_{\mathbb B}^{\frac 1 q}\Big({\hat\nu}[-M, x]^{1-p}+{\hat\nu}[y, N]^{1-p}\Big)^{-\frac 1 p}\bigg\}.
\endaligned$$
\end{itemize}
Moreover, we have $B_{\mathbb B*}\le B_{\mathbb B}^*\le 2^{1/p-1/q}B_{\mathbb B*}$ whenever $q\ge p$.
}
\dethm

\prf Let $g\in {\scr G}$. Without loss of generality, assume that $g>0$. For the
pair $\mu_g:=g \mu$ and $\nu$, by Theorem \ref{t1-7}, we know that the corresponding optimal constant $A_g$ in (\ref{01}) satisfies
$$B_{g*}\le A_g\le k_{q, p}\, B_g^*,$$
where
$$\aligned
&B_{g}^*=\sup_{-M\le x\le y\le N} \bigg\{\mu_g [x, y]^{\frac 1 q}\Big({\hat\nu}[-M, x]^{\frac{q(1-p)}{p}}+{\hat\nu}[y, N]^{\frac{q(1-p)}{p}}\Big)^{-\frac 1 q}\bigg\},\\
&B_{g*}=\sup_{-M\le x\le y\le N} \bigg\{\mu_g [x, y]^{\frac 1 q}\Big({\hat\nu}[-M, x]^{1-p}+{\hat\nu}[y, N]^{1-p}\Big)^{-\frac 1 p}\bigg\}.
\endaligned$$
Hence
$$\aligned
\sup_{g\in {\scr G}} B_{g*}
&= \sup_{x\le y} \sup_g \bigg\{\mu_g [x, y]^{1/q}\Big({\hat\nu}[-M, x]^{1-p}+{\hat\nu}[y, N]^{1-p}\Big)^{-1/p}\bigg\}\\
&= \sup_{x\le y} \bigg\{\Big(\sup_g \mu_g [x, y]\Big)^{1/q}\Big({\hat\nu}[-M, x]^{1-p}+{\hat\nu}[y, N]^{1-p}\Big)^{-1/p}\bigg\}\\
&= \sup_{x\le y} \bigg\{\big\|\mathbbold{1}_{ [x, y]}\big\|_{\mathbb B}^{1/q}\Big({\hat\nu}[-M, x]^{1-p}+{\hat\nu}[y, N]^{1-p}\Big)^{-1/p}\bigg\}\\
&= B_{\mathbb B*}.
\endaligned
$$
Similarly, we have
$\sup_{g\in {\scr G}} B_{g}^* = B_{\mathbb B}^*$.
From these facts, we obtain the estimates of $A_{\mathbb B}=\sup_{g\in {\scr G}}A_g$ immediately. The last assertion then follows from Lemma \ref{t1-4-1} (or its proof).
\deprf

As an application of Theorem \ref{t1-9}, it follows that the optimal constant $A_{\mathbb B}$ in the inequality
$$\|f^p\|_{\mathbb B}^{1/p}\le A_{\mathbb B}\|f'\|_{\nu, p},\qqd f(-M)=0\text{\cms\;and } f(N)=0$$
satisfies
\be B_{\mathbb B*}\le A_{\mathbb B}\le
k_{p, p}\, B_{\mathbb B*}, \lb{10-0}\de
where
$$B_{\mathbb B*}=\sup_{-M\le x\le y\le N} \bigg\{\big\|\mathbbold{1}_{ [x, y]}\big\|_{\mathbb B}^{\frac 1 p}\Big({\hat\nu}[-M, x]^{1-p}+{\hat\nu}[y, N]^{1-p}\Big)^{-\frac 1 p}\bigg\}.$$
Clearly, this result is simpler than Theorem \ref{t1-9}.
 Now, applying this result to ${\mathbb B}=L^{q/(q-p)}\,(q>p$, where $q/(q-p)$ is the conjugate number of
$q/p$), we return to Theorem \ref{t1-7} with a
factor $k_{p, p}$ different from $k_{q, p}$. In other words, we have arrived at
 a conclusion that the optimal constant
$A$ in the Hardy-type inequality (\ref{01}) (with $q\ge p$)
satisfies
\be B_*\le A\le k_{p, p}\, B_*,\lb{10-1}\de
where $B_*$ is the same as in Theorem \ref{t1-7}.
However, as we will see from Example \ref{t1-12} below that this result (\ref{10-1}) may be less sharp than Theorem \ref{t1-7}. Thus, lifting the left-hand side of the Hardy-type inequality from $L^p(\mu)$ to ${\mathbb B}$ does keep the constant $k_{p, p}$ but can not improve it.
We have thus explained the reason why on the left-hand side of the first inequality in
Theorem \ref{t1-9}, we use $q$ but not $p$.
In the discrete context and $p=2$, the conclusion (\ref{10-0}) was presented
by \rf{cmf10}{Theorem 8.2}.

In view of the direct proof for the lower estimate of $A$, the restriction $\mu_{\text{\rm pp}}=0$ coming from our splitting technique seems unnecessary. As just mentioned, the conclusion (\ref{10-0}) was proved in the discrete situation when $p=2$ by \rf{cmf10}{Theorem 8.2}. Its proof should be meaningful in the present continuous case (cf. \rf{cmf05}{Proof of Corollary 7.6}, simply replacing the function $e^{C(x)}$ there by the one $h(x)$ here). Hence for (\ref{10-0})
in the case of $p=2$, the condition $\mu_{\text{\rm pp}}=0$ is again not needed.
Thus, one may remove this condition by the capacitary method first
for $q=p$, if possible, and then extend to general normed linear space $\mathbb B$ as shown above. The present
proof depends heavily on the splitting property as shown by the first step of the proof of Lemma \ref{t1-3} and the comments below
Lemma \ref{t1-3}.
This is one of the reasons why the case of $q<p$ is missed here.
Actually, the last case is a rather different story, refer to
\rf{ok90}{Theorem 8.17} or \rf{mav}{pages 50-51}.

To illustrate the application of Theorem \ref{t1-7}, we study two examples.

\xmp\lb{t-11}{\cms Let $(-M, N)=(0, 1)$ and $\d\mu=\d\nu =\d x$. Then $H^{\partial}=0$ and
$$\aligned
B_*&=H^o= \frac 1 2 \bigg(\frac{p}{p - q + p q}\bigg)^{\frac 1 q}
\bigg(\frac{(p-1) q}{p - q + p q}\bigg)^{\frac{p-1}{p}},\\
B^*&=2^{1/p-1/q}H^o,\qqd 1<p\le q<\infty.
\endaligned$$
When $p=q=2$, it is known that $A=\pi^{-1}$ and $B^*=B_*=1/4$
(cf. \rf{cmf11}{Example 5.2}).}
\dexmp

\prf We have $h\equiv 1$, $y(x)=1-x$, and $m(\hat\nu)=1/2$. The
function
$$H_{\mu, \nu}(x, y(x))=2^{-\frac 1 p} \mu[x, y(x)]^{\frac 1 q}\,{\hat\nu}[0, x]^{\frac{p-1}{p}}=2^{-\frac 1 p} (1-2 x)^{\frac 1 q}x^{\frac{p-1}{p}}$$
achieves its maximum at
$$x=\frac 1 2\bigg(1+\frac{p}{(p-1)q}\bigg)^{-1}< \frac 1 2.$$
From this, we obtain $H^o$. Next, note that
$$\aligned
B_*&=\sup_{x\le y}\frac{(1-(1-y)-x)^{1/q}}{(x^{1-p}+(1-y)^{1-p})^{1/p}},\\
B^*&=\sup_{x\le y}\frac{(1-(1-y)-x)^{1/q}}{(x^{q(1-p)/p}+(1-y)^{q(1-p)/p})^{1/q}}.
\endaligned$$
Both of them are symmetric in $x$ and $1-y$. Hence $B_*=H^o$ and furthermore
$B^*=2^{1/p-1/q}H^o$.
\deprf

The next example illustrates the role played by $H^{\partial}$.

\xmp\lb{t1-12}{\cms Let $\mu(\d x)=\d x$ and $\nu(\d x)= x^2 \d x$ on $(1, \infty)$.
Assume that $p\in (1, 3)$. Then the inequality does not hold if $q\in [p, p/(3-p))$.
Otherwise, the inequality holds with
$$\aligned
&B_*=B^*=H^{\partial}=\bigg(\frac{p-1}{3-p}\bigg)^{\frac{p-1}{p}}\qqd\text{\cms if }q=\frac{p}{3-p}\text{\cms\; and }p\in [2, 3);\\
&\text{\rm $B_*$ and $B^*$ are bounded in terms of $H^o$ if $q>\frac{p}{3-p}$ and $q\ge p$},
\endaligned$$
where
$$H^o=2^{-\frac 1 p}\bigg(\frac{p-1}{3-p}\bigg)^{\frac{p-1}{p}}\!\!
\sup_{x\in (1, \, 2^{(3-p)/(p-1)}]}\bigg[\Big(1-x^{\frac{p-3}{p-1}}\Big)^{\frac{p-1}{p-3}}-x\bigg]^{\frac 1 q}\Big(1-x^{\frac{p-3}{p-1}}\Big)^{\frac{p-1}{p}}\!\!.$$
When $p=q=2$, we have $B^*=B_*=1$. In this case, $A=2$ and so the upper estimate $2B^*$ in Theorem \ref{t1-7} is exact (cf. \rf{cmf11}{Example 5.4}).
For fixed $p=2$, when $q$ varies from $2.01$ to $4.8$, the five
quantities we have worked so far are shown in Figure 1. The ratio of the upper and lower bounds is decreasing in $q$ but no more than $2$.
Note that we have a common lower bound $B_*$ in Theorem \ref{t1-7} and (\ref{10-1}),
but the upper bound in Theorem \ref{t1-7} is better than that in (\ref{10-1}).
\begin{center}{\includegraphics[width=11.0cm,height=7.5cm]{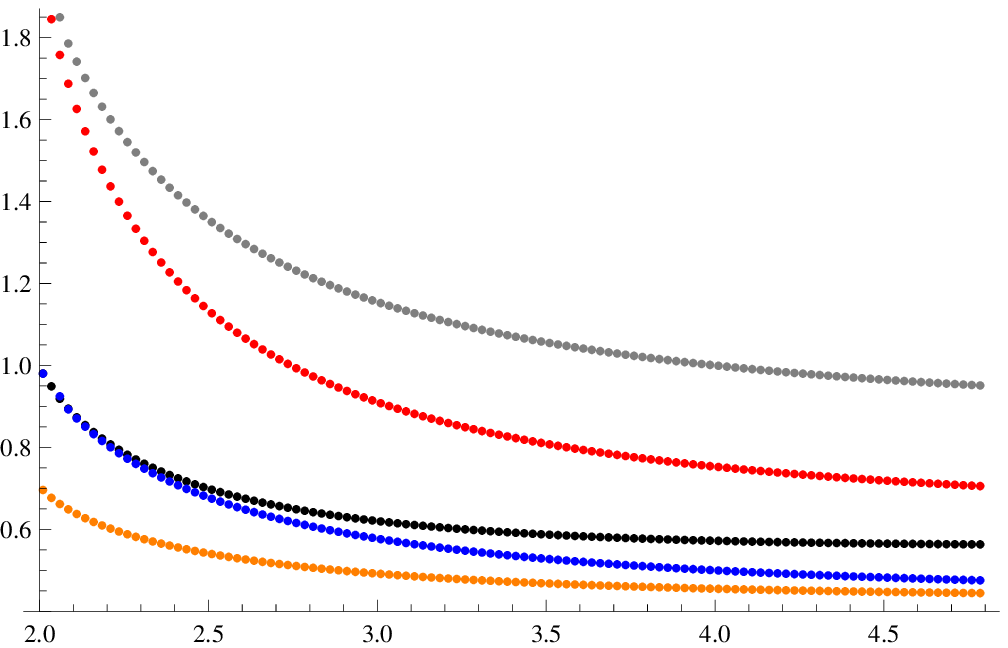}\newline
{\bf Figure 1}\qd The curves from bottom to top are
$H^o$, $B_*$, $B^*$, $k_{q, p} B^*$, and $k_{p, p}B_*$ respectively.}\end{center}
}
\dexmp

\prf We have $h(x)=x^{-2/(p-1)}$. Then
$$\int_1^x h=
\frac{p-1}{3-p}\Big(1-x^{\frac{p-3}{p-1}}\Big),\qqd
\int_y^{\infty} h=
\frac{p-1}{3-p}\, y^{\frac{p-3}{p-1}},
$$
where and in what follows, the Lebesgue measure $\d z$ is omitted.
Hence
\begin{gather}y(x)=\Big(1-x^{\frac{p-3}{p-1}}\Big)^{\frac{p-1}{p-3}},
\qqd m(\hat\nu)=2^{\frac{3-p}{p-1}},\nnb\\
H_{\mu, \nu}(x, y)=\bigg(\frac{p-1}{3-p}\bigg)^{\frac{p-1}{p}}
(y-x)^{\frac 1 q}\bigg[\Big(1-x^{\frac{p-3}{p-1}}\Big)^{1-p}
+y^{3-p}\bigg]^{-\frac 1 p}. \nnb
\end{gather}
In particular,
$$\varlimsup_{y\to\infty}H_{\mu, \nu}(x, y)
\!=\!\lim_{y\to\infty}\!\bigg(\frac{p-1}{3-p}\bigg)^{\frac{p-1}{p}} y^{1+\frac 1 q -\frac 3 p}
\!=\!\!
\begin{cases}
\Big(\frac{p-1}{3-p}\Big)^{\frac{p-1}{p}}&\text{ if $1+\frac 1 q = \frac 3 p$}\\
\infty &\text{ if $1+\frac 1 q > \frac 3 p$}\\
0 &\text{ if $1+\frac 1 q < \frac 3 p$}.
\end{cases}$$
The right-hand side is our $H^{\partial}$. Thus, if $1+\frac 1 q> \frac 3 p$, then
$B_*=\infty$. Next, we have
$$H_{\mu, \nu}(x, y(x))=2^{-\frac 1 p}\bigg(\frac{p-1}{3-p}\bigg)^{\frac{p-1}{p}}
\bigg[\Big(1-x^{\frac{p-3}{p-1}}\Big)^{\frac{p-1}{p-3}}-x\bigg]^{\frac 1 q}\Big(1-x^{\frac{p-3}{p-1}}\Big)^{\frac{p-1}{p}}.$$
Then
$$H^o=2^{-\frac 1 p}\bigg(\frac{p-1}{3-p}\bigg)^{\frac{p-1}{p}}
\sup_{x\in \big(1, \, 2^{(3-p)/(p-1)}\big]}\bigg[\Big(1-x^{\frac{p-3}{p-1}}\Big)^{\frac{p-1}{p-3}}-x\bigg]^{\frac 1 q}\Big(1-x^{\frac{p-3}{p-1}}\Big)^{\frac{p-1}{p}}.$$
The point here is that $H^o\le H^{\partial}$ and
$2^{1/p}H^o\le H^{\partial}$  in the case of
$q={p}/(3-p)$ and $p\ge 2$.
Besides, we have
$$\aligned
B_*&=\bigg(\frac{p-1}{3-p}\bigg)^{\frac{p-1}{p}}
\sup_{x<y}\left\{(y-x)^{\frac 1 q}\bigg[\Big(1-x^{\frac{p-3}{p-1}}\Big)^{1-p}
+y^{3-p}\bigg]^{-\frac 1 p}\right\},\\
B^*&=\bigg(\frac{p-1}{3-p}\bigg)^{\frac{p-1}{p}}
\sup_{x<y}\left\{(y-x)^{\frac 1 q}\bigg[\Big(1-x^{\frac{p-3}{p-1}}\Big)^{\frac{q(1-p)}{p}}
+y^{\frac{q(3-p)}{p}}\bigg]^{-\frac 1 q}\right\}.
\endaligned$$
Finally, numerical computation gives us the quantities $B^*$ et al, as shown in Figure 1.
\deprf

\section{The case of mean zero.}

Throughout this section, we assume that $\mu[-M, N]<\infty$ and
define a probability measure $\pi=(\mu [-M, N])^{-1}\mu$. In
probabilistic language, we are going to study the ergodic case
of the corresponding processes.
Corresponding to the three inequalities given at the beginning
of the last section, we now study the inequality
\be\bigg(\int_{-M}^N |f-\pi(f)|^q\d\mu\bigg)^{1/q}\le A
 \bigg(\int_{-M}^N |f'|^p\d\nu\bigg)^{1/p},\lb{11}\de
where $\pi(f)=\int f \d\pi$, in terms of
$$\aligned
&\bigg(\int_{\uz}^N |f|^q\d\mu\bigg)^{1/q}\le A_{\uz}^+ \bigg(\int_{\uz}^N |f'|^p\d\nu\bigg)^{1/p},\qquad f(\uz)=0,\\
&\bigg(\int_{-M}^{\uz} |f|^q\d\mu\bigg)^{1/q}\le A_{\uz}^- \bigg(\int_{-M}^{\uz} |f'|^p\d\nu\bigg)^{1/p},\qquad f(\uz)=0.
\endaligned$$
To save our notation, without any confusion, we use the same notation
$A$, $A_{\uz}^{\pm}$ and so on as in the last section.

Before moving further, let us mention the spectral meaning of
(\ref{01}) and (\ref{11}). Suppose that $\mu\ll\d x$ and
$\nu\ll\d x$, denote by $u=\d \mu/\d x$ and $v=\d \nu/\d x$.
Then the inverse of the optimal constant $A$ in (\ref{01}) and (\ref{11}), when $q=p$, corresponds to the infimum $\lz^{1/p}$ of the nontrivial spectrum of
$$\big(v |f'|^{p-1} \text{\rm sgn}(f')\big)'=-\lz u|f|^{q-1}\text{\rm sgn}(f)$$
with boundary condition $f(-M)=0=f(N)$ and $f'(-M)=0=f'(N)$
(when $M, N<\infty$),
respectively. The word ``bilateral'' in the title means that
a same boundary condition is endowed at two endpoints of the interval.
The spectral point of view has played a crucial
role in our previous study. For instance, it appears in each of
the papers \ct{cmf00} -- \ct{czz03}.

To study (\ref{11}), we start again at the splitting technique. We begin with the easier case: the lower estimate. It is indeed easier than the
one studied in the last section.

\lmm\lb{t2-1}{\cms Let $1\le p\le q< \infty$. Then we have
$$A\ge 2^{1/q-1/p} \sup_{\uz\in [-M, N]}\big(A_{\uz}^- \wedge A_{\uz}^+\big).$$}
\delmm

\prf Fix $\uz\in [-M, N]$. Let $f_-$ satisfy $f_-|_{[\uz,\, N]}=0$,
$$\int_{-M}^{\uz} |f_-|^q\d\mu=1,\qquad\text{and } \bigg(\int_{-M}^{\uz} |f_-'|^p\d\nu\bigg)^{1/p}< \big(A_{\uz}^-\big)^{-1}+\vz.$$
Let $f_+$ satisfy $f_+|_{[-M,\, \uz]}=0$,
$$\int_{\uz}^N |f_+|^q\d\mu=1,\qquad\text{and } \bigg(\int_{\uz}^N |f_+'|^p\d\nu\bigg)^{1/p}< \big(A_{\uz}^+\big)^{-1}+\vz.$$
Set $f=c f_- + f_+$, where $c=-\pi(f_+)/\pi(f_-)$. Then $\pi(f)=0$,
$$1+|c|^q=\int_{-M}^{\uz} |cf_-|^q\d\mu+\int_{\uz}^N |f_+|^q\d\mu=\int_{-M}^{N} |f|^q\d\mu,$$ and
$$\aligned
\int_{-M}^N |f'|^p\d\nu&=|c|^p\int_{-M}^{\uz} |f_-'|^p\d\nu+\int_{\uz}^N |f_+'|^p\d\nu\\
&\le |c|^p \Big(\big(A_{\uz}^-\big)^{-1}+\vz\Big)^p+ \Big(\big(A_{\uz}^+\big)^{-1}+\vz\Big)^p\\
&\le \Big(\big(A_{\uz}^-\big)^{-1}\vee \big(A_{\uz}^+\big)^{-1}+\vz\Big)^p (1+|c|^p).
\endaligned$$
Hence
$$\aligned
\bigg(\int_{-M}^N |f'|^p\d\nu\bigg)^{1/p}
&\le  \Big(\big(A_{\uz}^-\big)^{-1}\vee \big(A_{\uz}^+\big)^{-1}+\vz\Big) (1+|c|^p)^{1/p}\\
&\le 2^{1/p-1/q} \Big(\big(A_{\uz}^-\big)^{-1}\vee \big(A_{\uz}^+\big)^{-1}+\vz\Big) (1+|c|^q)^{1/q}\\
&\qquad\qquad\quad\text{(by Jensen's inequality requiring $q\ge p$)}\\
&=2^{1/p-1/q} \Big(\big(A_{\uz}^-\big)^{-1}\vee \big(A_{\uz}^+\big)^{-1}+\vz\Big) \bigg(\int_{-M}^{N} |f|^q\d\mu\bigg)^{1/q}.
\endaligned$$
Thus
$$A\ge 2^{1/q-1/p} \big(A_{\uz}^- \wedge A_{\uz}^+\big).$$
Since $\uz$ is arbitrary, we obtain the lower bound of $A$.
\deprf

The upper bound of $A$ is harder
than the lower one just studied.
But the first step is still easy. Given $f$ and $\uz\in (-M, N)$, let $\tilde f=f-f(\uz)$. Then
$$\aligned
\int_{-M}^N \big|{f}'\big|^p\d\nu&=\int_{-M}^{\uz} \big|{\tilde f}'\big|^p\d\nu+\int_{\uz}^N \big|{\tilde f}'\big|^p\d\nu\\
&\ge  \big(A_{\uz}^-\big)^{-p}\bigg(\int_{-M}^{\uz} \big|\tilde f\big|^q\d\mu\bigg)^{p/q}+  \big(A_{\uz}^+\big)^{-p}\bigg(\int_{\uz}^N \big|\tilde f\big|^q\d\mu\bigg)^{p/q}\\
&\ge \Big[\big(A_{\uz}^-\big)^{-p}\wedge \big(A_{\uz}^+\big)^{-p}\Big]\bigg[\bigg(\int_{-M}^{\uz} \big|\tilde f\big|^q\d\mu\bigg)^{p/q}
  + \bigg(\int_{\uz}^N \big|\tilde f\big|^q\d\mu\bigg)^{p/q}\bigg]\\
&\ge \big(2^{(p/q-1)\vee 0}\big)^{-1}\Big[\big(A_{\uz}^-\big)^{-p}\wedge \big(A_{\uz}^+\big)^{-p}\Big] \bigg(\int_{-M}^{N} \big|\tilde f\big|^q\d\mu\bigg)^{p/q}\\
&\qquad\qquad\quad(\text{by $c_r$-inequality}).
\endaligned$$
Our aim is to replace $\big|\tilde f\big|^q$ on the right-hand side with $|f-\pi(f)|^q$. This is true in the case of $q=2$ since
$$\inf_{c\in {\mathbb R}} \int_{-M}^{N} (f-c)^2\d\mu=\int_{-M}^{N} (f-\pi(f))^2\d\mu.$$
Unfortunately, it does not work for general $q$. Anyhow, when $q=2$, we have
$$\int_{-M}^N \big|{f}'\big|^p\d\nu \ge \big(2^{(p/2-1)\vee 0}\big)^{-1}\Big[\big(A_{\uz}^-\big)^{-p}\wedge \big(A_{\uz}^+\big)^{-p}\Big] \bigg(\int_{-M}^{N} \big|f-\pi(f)\big|^2\d\mu\bigg)^{p/2}.$$
Since $\uz$ is arbitrary, we obtain
$$\aligned
\big(2^{(1/2-1/p)\vee 0}\big)\bigg(\int_{-M}^N \big|{f}'\big|^p\d\nu\bigg)^{1/p} \ge &
\sup_{\uz\in (-M, N)}
\Big[\big(A_{\uz}^-\big)^{-1}\wedge \big(A_{\uz}^+\big)^{-1}\Big]\\
&\times \bigg(\int_{-M}^{N} \big|f-\pi(f)\big|^2\d\mu\bigg)^{1/2}.
\endaligned$$
Next, since $f$ is arbitrary, it follows that
$$A\le 2^{(1/2-1/p)\vee 0} \inf_{\uz\in (-M,\, N)}\big(A_{\uz}^- \vee A_{\uz}^+\big).$$

Up to now, the proof is similar to \rf{cmf00}{Theorems 3.3 and 3.4} in the specific case that $q=2$.
For general $q\ge 2$, we have luckily a different approach (cf. \rf{cmf05}{Chapter 6} and references therein). Note that we have already proved that
if the measure $\mu$ is replaced by $\mu_g:=g \mu$ for a nonnegative function $g$ on $[-M, N]$,
then the optimal constant $A_g$ in the inequality
$$\bigg(\int_{-M}^{N} \big|f-\pi(f)\big|^2\d\mu_g\bigg)^{1/2}\le A_g \bigg(\int_{-M}^N \big|{f}'\big|^p\d\nu\bigg)^{1/p}$$
obeys
$$A_g\le 2^{(1/2-1/p)\vee 0} \inf_{\uz\in (-M,\, N)}\big(A_g^{\uz, -} \vee A_g^{\uz, +}\big),$$
where $A_g^{\uz, \pm}$ is obtained from $A_{\uz}^{\pm}$ replacing $\mu$ with $\mu_g$.
From now on in this proof, the constants $A_g$, $A_g^{\uz, \pm}$, and $B_g^{\uz, \pm}$ are used for
$\mu_g$ in the specified case that $q=2$ only. Note that for $A_g^{\uz, +}$ for instance, the function $g$
can be replaced by $g\mathbbold{1}_{[\uz, N]}$.

Even though we are now mainly working on the $L^q$-case to which ${\scr G}$ is the set of functions in the unit ball of $L^{\frac{q}{q-2}} (\mu)$
\big(where $\frac{q}{q-2}$ is the conjugate number of $q/2$\big):
$$\aligned
\bigg(\!\int_{-M}^{N} \!\big|f\!-\!\pi(f)\big|^q\d\mu\bigg)^{2/q}\!\!\!&=\sup_{g\in {\scr G}} \int_{-M}^N \!\big|f\!-\!\pi(f)\big|^2 g\d\mu
\!=\!\sup_{g\in {\scr G}} \int_{-M}^N \!\big|f\!-\!\pi(f)\big|^2 \d\mu_g,
\endaligned$$
at the moment, we allow ${\scr G}$ to be general in the setup of Hypotheses \ref{t1-8}:
\be \big\|\big(f-\pi(f)\big)^2\big\|_{\mathbb B}^{1/2}\le A_{\mathbb B} \|f'\|_{\nu, p}.\lb{12}\de
We have
$$\aligned
A_{\mathbb B}&=\sup_{g\in {\scr G}}A_g\le  \sup_{g\in {\scr G}}\,\inf_{\uz\in (-M, N)}\big(A_g^{\uz, -} \vee A_g^{\uz, +}\big)\\
&\le \inf_{\uz\in (-M, N)}\bigg[\bigg(\sup_{g\in {\scr G}_{\uz}^-} A_g^{\uz, -}\bigg)
 \bigvee\bigg(\sup_{g\in {\scr G}_{\uz}^+} A_g^{\uz, +}\bigg)\bigg],\endaligned$$
where
$${\scr G}_{\uz}^{-}=\{g|_{[-M, \uz]}: g\in {\scr G}\},\qqd
{\scr G}_{\uz}^{+}=\{g|_{[\uz, N]}: g\in {\scr G}\}.$$

Here is a technical point.
Because on the left-hand side of (\ref{12}), we start at $q=2$. This leads to
the restriction that $p\in (1, 2]$ since we need $q\ge p$ in order
to use the basic estimates in terms of $B_{\uz}^{\pm}$ given in the next proof. Anyhow, we have proved the first assertion of the next result.

\lmm\lb{t2-2}{\cms Let ${\scr G}$ satisfy Hypotheses \ref{t1-8}. Then for $p\in (1, 2]$, we have
$$A_{\mathbb B} \le \inf_{\uz\in (-M, N)}\bigg[\bigg(\sup_{g\in {\scr G}_{\uz}^-} A_g^{\uz, -}\bigg)
 \bigvee\bigg(\sup_{g\in {\scr G}_{\uz}^+} A_g^{\uz, +}\bigg)\bigg].$$
Moreover,
$$\sup_{g\in {\scr G}_{\uz}^{\pm}}A_g^{\uz, \pm}\le k_{2, p}\, B_{\mathbb B}^{\uz, \pm},$$
where
$$\aligned
B_{\mathbb B}^{\uz, +}&= \sup_{r\in (\uz,\, N)} \big\|\mathbbold{1}_{[r,\, N]}\big\|_{\mathbb B}^{1/2} \,{\hat\nu}[\uz, r]^{(p-1)/p},\\
B_{\mathbb B}^{\uz, -}&=\sup_{r\in (-M,\, \uz)} \big\|\mathbbold{1}_{[-M,\, r]}\big\|_{\mathbb B}^{1/2} \,{\hat\nu}[r, \uz]^{(p-1)/p}.
\endaligned$$
 }\delmm

\prf By \rf{ok90}{Theorem 1.14} and \rf{mav}{\S 1.3, Theorem 1} (see also \rf{bg91}{Theorem 8}  and \rf{mvm92}{Theorem 2} in which the factor $ k_{q, p}$ may be different), we have for general $1<p\le q<\infty$ that
$$\aligned
B_{\uz}^{\pm}&\le A_{\uz}^{\pm}\le k_{q, p}\,B_{\uz}^{\pm},\qqd
1<p\le q<\infty,\\
B_{\uz}^+&= \sup_{r\in (\uz,\, N)} \mu [r, N]^{1/q}\,{\hat\nu}[\uz, r]^{(p-1)/p},\\
B_{\uz}^-&=\sup_{r\in (-M, \,\uz)} \mu [-M, r]^{1/q}\,{\hat\nu}[r, \uz]^{(p-1)/p},
\endaligned$$
where $\hat\nu$ is the same as in the last section.
It remains to estimate $\sup_{g\in {\scr G}_{\uz}^+} A_g^{\uz, +}$ for instance. First, we have
for $q=2$ that
$$\sup_{g\in {\scr G}_{\uz}^+} A_g^{\uz, +}
\le k_{2, p}\,\sup_{g\in {\scr G}_{\uz}^+} B_g^{\uz, +}.$$
Next, we have
$$\aligned
\sup_{g\in {\scr G}_{\uz}^+}B_g^{\uz, +}&=\sup_{g\in {\scr G}_{\uz}^+}\;\sup_{r\in (\uz,\, N)}\mu_g [r, N]^{1/2}\,{\hat\nu}[\uz, r]^{(p-1)/p}\\
&=\sup_{r\in (\uz,\, N)}\Big(\sup_{g\in {\scr G}_{\uz}^+}\mu_g [r, N]\Big)^{1/2}\,{\hat\nu}[\uz, r]^{(p-1)/p}\\
&= B_{\mathbb B}^{\uz, +}.\endaligned$$
Similar computation holds for $\sup_{g\in {\scr G}_{\uz}^-}B_g^{\uz, -}$. Combining these facts with the first
assertion gives us the second one of the lemma.
\deprf

As mentioned in the last section, without loss of generality, we
can assume that $\mu$ is positive on each subinterval.

Here is our upper estimate.

\lmm\lb{t2-3}{\cms Let $\mu_{\text{\rm pp}}=0$ and ${\scr G}$ satisfy Hypotheses \ref{t1-8}. Then for $p\in (1, 2]$, we have
$$
A_{\mathbb B} \le k_{2, p}\,{B_{\mathbb B}^*},$$
where the constant $B_{\mathbb B}^*$ is defined by
$${B_{\mathbb B}^*}^{-\frac{p}{p-1}}=\inf_{x< y}\Big[
\|\mathbbold{1}_{[-M,\, x]}\|_{\mathbb B}^{-\frac{p}{2(p-1)}}
+\|\mathbbold{1}_{[y,\, N]}\|_{\mathbb B}^{-\frac{p}{2(p-1)}} \Big]
\,{\hat\nu}[x, y]^{-1}.$$
}\delmm

\prf Write
$${B_{\mathbb B}^*}^{\frac{p}{p-1}}
=\sup_{x< y}
\frac{{\hat\nu}[x, \uz] +{\hat\nu}(\uz, y]}
{
\|\mathbbold{1}_{[-M,\, x]}\|_{\mathbb B}^{-\frac{p}{2(p-1)}}
+\|\mathbbold{1}_{[y,\, N]}\|_{\mathbb B}^{-\frac{p}{2(p-1)}} }.$$
Similar to the proof of Lemma \ref{t1-3}, the assertion follows by using Lemmas \ref{t1-2} and \ref{t2-2}. Here we may
need the approximating procedure by finite $M$ and $N$.
\deprf

The following result is on the lower estimate of $A$. Its proof is new even in the special case that $p=q=2$. Note that $\mu[x, y]=\mu (x, y)$ whenever $\mu_{\text{\rm pp}}=0$.

\lmm\lb{t2-4} {\cms Let $\mu_{\text{\rm pp}}=0$. Then for $1<p,\, q< \infty$, the optimal constant $A$ in (\ref{11}) satisfies
\begin{align}
A&\ge \!\!\sup_{-M\le x <y\le N}\bigg\{\Big[\mu[-M, x]^{\frac 1 {1-q}}+ \mu[y, N]^{\frac 1 {1-q}}\Big]^{\frac {1-q} q}{\hat\nu}[x, y]^{\frac{p-1}{p}}\bigg\}\nnb\\
&=:B_*.\lb{13}\end{align}
}
\delmm

\prf Given $m, n\in (-M,\,N)$ with $m<n$, let $\bar\uz=\bar\uz (m, n)$ be the unique solution to the equation
$$\mu[-M, m]\,{\hat\nu}[m, \uz]+ \!\int_m^{\uz}\! \!\!\mu(\d x){\hat\nu}[x, \uz]
\!=\!\mu[n, N]\,{\hat\nu}[\uz, n]+\!\int_{\uz}^n\!\!\! \mu(\d x){\hat\nu}[\uz, x],\;\; \uz\in (m, n).$$
The existence of the solution is clear since $\mu$ is continuous,
when $\uz$ varies from $m$ to $n$,
the left-hand side goes from $0$ to a positive number and the right-hand side goes from a positive number to zero.
Next, define
$$f(x)=- \mathbbold{1}_{\{x\le \bar\uz\}}{\hat\nu}\big[m\vee x, \bar\uz\,\big]
+ \mathbbold{1}_{\{x> \bar\uz\}} {\hat\nu}\big[\bar\uz, n\wedge x\big],\qqd x\in [-M, N].$$
Then $\mu(f)= 0$ by definition of $\bar\uz$.
Clearly, $f$ is absolutely continuous. On the one hand, we have
\begin{align}
\bigg(\int_{-M}^N |f'|^p\d\nu\bigg)^{1/p} &= \bigg( \int_m^{\bar\uz} h^p\d\nu +\int_{\bar\uz}^n h^p \d \nu\bigg)^{1/p}\nnb\\
&=\big({\hat\nu}\big[m, \bar\uz\,\big]  +{\hat\nu}\big[\bar\uz, n\big]\big)^{1/p}\nnb\\
&={\hat\nu}[m, n]^{1/p}.\lb{14}
\end{align}
Here in the second step, we have once again ignored the singular part of $\nu$. On the other hand, we have
\begin{align}
\int_{-M}^N |f-\pi (f)|^q\d\mu
&=\int_{-M}^N |f|^q\d\mu\nnb \\
&> \int_{-M}^m |f|^q\d\mu +\int_{n}^N |f|^q\d\mu\nnb\\
&=\mu[-M, m]\, {\hat\nu}\big[m, \bar\uz\,\big]^q
+\mu[n, N]\,{\hat\nu}\big[\bar\uz, n\big]^q.\lb{14-1}
\end{align}
Now we have naturally, as in proof (a) of Proposition \ref{t1-1}, that
$$\aligned
\text{RHS of (\ref{14-1})}&\ge \big(\mu[-M, m]\wedge \mu[n, N]\big)
\big({\hat\nu}\big[m, \bar\uz\,\big]^q+{\hat\nu}\big[\bar\uz, n\big]^q\big)\\
&\ge 2^{1-q}\big(\mu[-M, m]\wedge \mu[n, N]\big){\hat\nu}\big[m, n\big]^q.
\endaligned$$
However, such a lower bound is quite rough for our purpose so we need a different approach.
Note that the function
$$\gz(x)=\az x^q + \bz (1-x)^q,\qqd x\in (0, 1),
\;\az>0, \; \bz>0,\; q\in (1, \infty)$$
achieves its minimum
$$\Big(\az^{\frac 1 {1-q}}+ \bz^{\frac 1 {1-q}}\big)^{1-q}
\qqd (\text{resp., $ \az\wedge \bz$ in the case of $q=1$})$$ at
$$x^*=\bigg[1+\bigg(\frac{\az}{\bz}\bigg)^{\frac 1 {q-1}}\bigg]^{-1}=\bz^{\frac 1{q-1}}\Big[\az^{\frac 1{q-1}}+\bz^{\frac 1{q-1}}\Big]^{-1}\in (0, 1).$$
Applying this result with
$$\az=\mu[-M, m],\qd \bz=\mu[n, N],\qd x={\hat\nu}\big[m, \bar\uz\,\big]
\big/{\hat\nu}[m, n] $$
to (\ref{14-1}), we get
$$\bigg(\int_{-M}^N |f-\pi (f)|^q\d\mu\bigg)^{1/q}
\ge\Big\{\mu[-M, m]^{\frac 1 {1-q}}+ \mu[n, N]^{\frac 1 {1-q}}\Big\}^{\frac {1-q} q}
{\hat\nu}[m, n]. $$
Because
$$A\ge \bigg(\int_{-M}^N |f-\pi (f)|^q\d\mu\bigg)^{1/q}\bigg(\int_{-M}^N |f'|^p\d\nu\bigg)^{-1/p}, $$
the estimate given in the lemma now follows immediately.
\deprf

Now, one may ask the possibility using the idea in the last part
of the proof above to improve the estimate produced by proof (a) of Proposition \ref{t1-1}. The answer is yes if $p>q$ and no if $p\le q$. Note that here we have power $q>1$ and in proof (a) of Proposition \ref{t1-1}, the power is $p/q$. Thus, if $p>q$, we can follow the proof here to have an improvement. However, in this paper, we are mainly interested in the case that $p\le q$. Then the function $\az x^{\gz}+ \bz (1-x)^{\gz}\,(\gz\le 1)$ is concave, its minimum is achieved at the boundaries: either at $x=0$ or at $x=1$. That is, $\min_{x\in (0, 1)}\{\az x^{\gz}+ \bz (1-x)^{\gz}\}=\az\wedge \bz$. In this case, we have thus returned to the original result given in proof (a) of Proposition \ref{t1-1}.
We mention that this remark is also meaningful for the first step of the proof of Lemma \ref{t2-2} given
right below the proof of Lemma \ref{t2-1}.

As an analog of Lemma \ref{t1-4-1}, we have the following result.

\lmm\lb{t2-4-1}{\cms Let $q\ge p$. Then we have $B_*\le B^*\le 2^{1/p-1/q}B_*$,
where the constant $B^*$ is defined by
\be {B^*}^{\frac{p}{1-p}}\!=\!\inf_{x< y}\Big\{
\mu [-M, x]^{\frac{p}{(1-p)q}}
\!\!+\!\mu [y, N]^{\frac{p}{(1-p)q}} \Big\}\,
{\hat\nu}[x, y]^{-1}\!\!.\lb{17}\de}
\delmm

\prf Applying ${\mathbb B}$ to $L^{\frac{q}{q-2}} (\mu)$, the constant $B_{\mathbb B}^*$
given in Lemma \ref{t2-3} is reduced to $B^*$ defined by (\ref{17}).

(a) Part $B^*\ge B_*$ follows from the $c_r$-inequality by setting
$$\az=\mu[-M, x]^{\frac{p}{(1-p)q}},\qd
\bz=\mu[y ,N]^{\frac{p}{(1-p)q}},\qd
r=\frac {(p-1) q}{p(q-1)} \in (0, 1].$$

(b) Part $B_*\ge 2^{1/q-1/p}B^*$ follows from the inequality by setting
$$\az=\mu[-M, x]^{\frac{1}{1-q}},\qd
\bz=\mu[y ,N]^{\frac{1}{1-q}},\qd
r=\frac {p(q-1)}{(p-1) q} \ge 1.\qqd\square$$

As a combination of Lemmas \ref{t2-3} -- \ref{t2-4-1}, we obtain the following result.

\thm\lb{t2-6}{\cms Let $\mu[-M, N]<\infty$ and $\mu_{\text{\rm pp}}=0$. Then
\begin{itemize}\setlength{\itemsep}{-0.8ex}
\item[(1)] for $1< p\le 2 \le q < \infty$, the optimal constant $A$ in (\ref{11}) satisfies
$$ A \le k_{2, p}\, {B^*},\;\text{\cms where $k_{q, p}$ is defined by (\ref{04-1}), and}$$
\item[(2)] for $1<p,\, q<\infty$, we have $A\ge B_*$,
\end{itemize}
where $B_*$ and $B^*$ are defined in Lemmas \ref{t2-4} and (\ref{17}), respectively. Moreover, we have $B_*\le B^*\le 2^{1/p-1/q}B_*$ once $q\ge p$.
}\dethm

\prf As an application of Lemma \ref{t2-3}, we
get the upper estimate of $A$.
The lower estimate of $A$ is due to Lemma \ref{t2-4}.
The comparison of $B^*$ and $B_*$ comes from Lemma \ref{t2-4-1}.
\deprf

When $p=q=2$, from Theorem \ref{t2-6}, it follows that
$$ B^*\le A \le 2 B^*.$$
We have thus returned to \rf{cmf10}{Theorem 10.2}. It is interesting that in the special case of $p=q=2$, the duality given in \rf{kp03}{page 13 and (1.17)} coincides with that used in \ct{cmf10, cmf11}. The former duality exchanges the (single-side but not bilateral) boundary conditions $f(-M)=0$ and $f(N)=0$. This is clearly different from a dual of (\ref{01}) and (\ref{11}).
To prove the last duality, in \ct{cmf10, cmf11}, several techniques were adopted: coupling, duality, and capacity. Thus,
the proofs given here are essentially different from that presented in \ct{cmf10, cmf11}, much direct and elementary. Besides, it is unclear how these advanced techniques can be applied to the present setup.

 In parallel to the last section, define $y(x)$ to be the solution to the equation $\mu[-M, x]=\mu[y, N]$ and denote by $m(\mu)$ to
 be the median of $\mu$. Set
 $$H_{\mu, \nu}(x, y)=\Big[
\mu[-M, x]^{\frac 1{1-q}}\!+\!\mu[y, N]^{\frac 1{1-q}} \Big]^{\frac {1-q}q}\,
{\hat\nu}[x, y]^{\frac{p-1}{p}}.$$
Define
\be H^o\!=2^{1/q-1}\!\!\sup_{x\in (-M,\, m(\mu)]}\mu[-M, x]^{1/q}\,
{\hat\nu}[x, y(x)]^{(p-1)/p}\!\!.\lb{15}\de
Denote by $\ggz$ be the limiting points of $H_{\mu, \nu}(x, y)$ as
$${\hat\nu}[-M, x]=\infty\qd\text{or}\qd
{\hat\nu}[y, N]=\infty,$$ as well as the iterated limits if
${\hat\nu}[-M, N]=\infty$
when $M=\infty=N$. Set
\be H^{\partial}=
\begin{cases}
\sup \{\gz: \gz\in \ggz\}\qd &\text{if } \ggz\ne\emptyset\\
0 \qd &\text{if } \ggz=\emptyset.
\end{cases}\lb{16}\de

Similar to Lemmas \ref{t1-5} and \ref{t1-6}, we have the following result.

\lmm\lb{t2-5}{\cms Let $\mu[-M, N]<\infty$ and $\mu_{\text{\rm pp}}\!=\!0$. Define $H^o$ and $H^{\partial}$ as above. Then we have
$$H^o \vee H^{\partial}\le B_*\le  \big(2^{1-1/q}H^o\big)\vee H^{\partial}$$ and
$$\big(2^{1/p-1/q}H^o\big)\vee H^{\partial}\le B^*\le \big(2^{1-1/q}H^o\big)\vee H^{\partial}.$$
}
\delmm

\prf Consider $B_*$ for instance. Recalling that $\mu [-M, x]= \mu [y(x), N]$, we have
$$\aligned
&\sup_{-M\le x<y\le N}\Big\{\mu[-M, x]^{\frac 1 {1-q}}+ \mu[y, N]^{\frac 1 {1-q}}\Big\}^{\frac {1-q} q}\,{\hat\nu}[x, y]^{\frac{p-1}{p}}\\
&\qd\ge2^{1/q-1} \sup_{-M\le x\le m(\mu)} \mu[-M, x]^{1/q}\, {\hat\nu}[x, y(x)]^{(p-1)/p}\\
&\qd=  H^{o}.
\endaligned$$
This plus the boundary condition gives us the lower estimate of $B_*$. The proofs for the other assertions are similar.
\deprf

\crl\lb{t2-6-1}{\cms The Hardy-type inequality (\ref{11}) holds iff $H^o\vee H^{\partial}<\infty$.}
\decrl

To generalize Theorem \ref{t2-6} to a general normed linear space ${\mathbb B}$, as in the study of
the upper estimate $B_{\mathbb B}^*$, a natural way is starting from $B_{g*}$:
$$\sup_g B_{g*}=\sup_g \sup_{x<y}\bigg\{\Big(\mu_g[-M, x]^{\frac 1{1-q}}+\mu_g[y, N]^{\frac 1{1-q}}\Big)^{\frac {1-q}q}\,{\hat\nu}[x, y]^{\frac{p-1}{p}}\bigg\}.$$
Then it is not clear how to handle with this expression
in terms of the norm $\mathbb B$. A crucial point here is that the measure $\mu$ appears in the last expression twice
rather than a single term in the last section.
The next result is an extension and improvement
of the basic estimates given in \rf{cmf02}{Theorem 2.2}.

\thm\lb{t2-9}{\cms Let $\mu[-M, N]<\infty$, $\mu_{\text{\rm pp}}=0$, and ${\scr G}$ satisfy Hypotheses \ref{t1-8}. Then
\begin{itemize}\setlength{\itemsep}{-0.8ex}
\item[(1)] for $p\in (1, 2]$, the optimal constant $A_{\mathbb B}$ in the inequality
$$ \|\big(f-\pi(f)\big)^2\|_{\mathbb B}^{1/2}\le A_{\mathbb B} \|f'\|_{\nu, p}$$
satisfies $A_{\mathbb B} \le k_{2, p}\,{B_{\mathbb B}^*},$
where the constant $B_{\mathbb B}^*$ is defined by
$${B_{\mathbb B}^*}^{-\frac{p}{p-1}}=\inf_{-M<x< y<N}\Big\{
\|\mathbbold{1}_{[-M, x]}\|_{\mathbb B}^{-\frac{p}{2(p-1)}}
+\|\mathbbold{1}_{[y, N]}\|_{\mathbb B}^{-\frac{p}{2(p-1)}} \Big\}
\,{\hat\nu}[x, y]^{-1}.$$
\item[(2)] For $1<p, q <\infty$, the
optimal constant $A_{\mathbb B}$ in the inequality
$$ \big\||f-\pi(f)|^q\big\|_{\mathbb B}^{1/q}\le A_{\mathbb B} \|f'\|_{\nu, p}$$
satisfies $A_{\mathbb B} \ge B_{\mathbb B *},$
where
$$\!B_{\mathbb B *}=\sup_{-M<x< y<N}\gz_{\mathbb B}^{}(x, y; q)\, {\hat\nu}[x, y]^{\frac{p-1}{p}}\!\!$$
and  $$\gz_{\mathbb B}^{}(x, y; q)=\inf_{z\in (0, 1)}
  \big\|\mathbbold{1}_{[-M, x]}z^q
  +\mathbbold{1}_{[y, N]}(1-z)^q\big\|_{\mathbb B}^{1/q}.$$
\end{itemize}
}\dethm

\prf The first assertion is a copy of Lemma \ref{t2-3}. To prove the second assertion,
we return to the construction used in the proof of Lemma \ref{t2-4}. That is, we use the
notation $\bar\uz$ and $f$ introduced there. First, we have
$$\aligned
\big\||f-\pi(f)|^q\big\|_{\mathbb B}
&= \big\||f|^q\big\|_{\mathbb B}\\
&\ge \big\| |f|^q \mathbbold{1}_{[-M, m]}+ |f|^q\mathbbold{1}_{[n, N]}\big\|_{\mathbb B}\\
&=\Big\| \mathbbold{1}_{[-M, m]}\,{\hat\nu}\big[m, \bar\uz\big]^q
+ \mathbbold{1}_{[n, N]}\,{\hat\nu}\big[\bar\uz, n\big]^q\Big\|_{\mathbb B}\\
&\ge \gz_{\mathbb B}^{}(m, n; q)^q \,{\hat\nu}[m, n]^q.\endaligned$$
Combining this with (\ref{14}), we obtain
$$\aligned
\frac{\big\||f\!-\!\pi(f)|^q\big\|_{\mathbb B}^{1/q}}{\|f'\|_{\nu, p}}
&\ge \gz_{\mathbb B}^{}(m, n; q)
 \, {\hat\nu}[m, n]^{\frac{p-1}{p}}\!\!.\endaligned$$
Now the required assertion follows by making supremum with respect to $(x, y)$ with $x<y$.
\deprf

\xmp\lb{t2-7}{\cms Let $\mu(\d x)=\nu(\d x)=e^{-b x}\d x$ ($b>0$) on $(0, \infty)$.
Then the inequality (\ref{11}) does not hold if $q>p$. When $q=p$, we have
$$B^*=B_*=H^{\partial}=\frac{1}{b}\big(p-1\big)^{1-\frac 1 p}.$$
The upper estimate $2B^*$ in Theorem \ref{t2-6}\,(1) is sharp in the case of $p=q=2$, refer to \rf{cmf11}{Example 5.3}.}
\dexmp

\prf We have
$$\mu(0, x)\!=\!\frac 1 b (1- e^{-b x}),\;
\mu(y, \infty)\!=\!\frac 1 b e^{-b y},\; y(x)\!=\!-\frac 1 b \log (1- e^{-b x}),\,
m(\mu)\!=\!\frac 1 b \log 2.$$
Next,
$$h(x)=e^{\frac{bx}{p-1}},\qqd
{\hat\nu}[x, y]= \frac{p-1}{b}\Big(e^{\frac{b y}{p-1}}- e^{\frac{b x}{p-1}}\Big).$$
Thus,
$$H_{\mu, \nu}(x, y)\!=\!b^{-\frac 1 q}\bigg(\frac{p-1}{b}\bigg)^{1-\frac 1 p}
\Big[(1- e^{-b x})^{\frac 1{1-q}}+e^{-\frac{b y}{1-q}}\Big]^{\frac 1 q-1}\Big(e^{\frac{b y}{p-1}}- e^{\frac{b x}{p-1}}\Big)^{1-\frac 1 p}.$$
Hence
$$\aligned
H^{\partial}\!&=\!\varlimsup_{y\to\infty}H_{\mu, \nu}(x, y)\!=\!b^{-\frac 1 q}\bigg(\frac{p-1}{b}\bigg)^{1-\frac 1 p}
\lim_{y\to\infty}e^{by \big(\frac 1 p - \frac 1 q\big)}\\
\!&=\!
{\begin{cases}
b^{-1} (p-1)^{1-\frac 1 p}&\text{if } q\!=\!p\\
\infty   &\text{if } q\!>\!p.
\end{cases}}
\endaligned$$
Thus, the inequality (\ref{11}) does not hold if $q>p$. Next, assume that $q=p$. Then we have
$$\aligned
H^o&\!=\!2^{\frac 1 p-1}b^{-\frac 1 p}\bigg(\frac{p-1}{b}\bigg)^{1-\frac 1 p}\!\!\!\!\!\sup_{x\in (0,\, b^{-1}\log 2)}\!(1- e^{-b x})^{\frac 1 p}\Big((1- e^{-b x})^{-\frac {1}{p-1}}\!-e^{\frac{b x}{p-1}}\Big)^{1-\frac 1 p}\\
&\!=\!2^{\frac 1 p-1}b^{-\frac 1 p}\bigg(\frac{p-1}{b}\bigg)^{1-\frac 1 p}\!\!\sup_{x\in (0,\, \log 2)}(1- e^{-x})^{\frac 1 p}\Big((1- e^{-x})^{-\frac {1}{p-1}}-e^{\frac{x}{p-1}}\Big)^{1-\frac 1 p}.
\endaligned$$
To compute $H^o$, we observe that
$$\varlimsup_{x\to 0}(1- e^{-x})^{\frac 1 p}\Big((1- e^{-x})^{-\frac {1}{p-1}}-e^{\frac{x}{p-1}}\Big)^{1-\frac 1 p}
=1.$$
Because of this and the decreasing property of the function on the left-hand side, we obtain
 $$H^o=2^{\frac 1 p-1}b^{-\frac 1 p}\Big(\frac{p-1}{b}\Big)^{1-\frac 1 p}=2^{\frac 1 p-1}
 \frac 1 b (p-1)^{1-\frac 1 p}.$$
Having $H^{\partial}$ and $H^{o}$ at hand, it is easy to compute $B^*$ and $B_*$.
\deprf

\xmp\lb{t2-8}{\cms Let $\mu(\d x)=x^{-2}\d x$ and $\nu(\d x)=\d x$ on $(1, \infty)$.
Then the inequality (\ref{11}) does not hold if $\frac 1 p + \frac 1 q<1$. Otherwise, if
$\frac 1 p + \frac 1 q=1$, then $B^*=B_*=H^{\partial}=1$.
In particular, when $p=q=2$, our upper
estimate $2B^*$ in Theorem \ref{t2-6}\,(1) is exact since $A=2$ (cf. \rf{cmf11}{Example 5.4}).
If $\frac 1 p + \frac 1 q>1$, then ${H}^o\le B_*\le 2^{1-1/q} {H}^o$ and
$2^{1/p-1/q}{H}^o\le B^*\le 2^{1-1/q}{H}^o$:
$${H}^o\!=\!2^{\frac 1 q-1}\!\left[\frac{2 \bz}{\sqrt{\az^2-6 \az \bz+\bz^2}+\az-\bz}+\!1\right]^\az \!\left[2-\frac{\sqrt{\az^2-6 \az
   \bz+\bz^2}+\az+\bz}{2 \bz}\right]^\bz\!\!\!,$$
   where
   $\az= 1-\frac 1 p - \frac 1 q<0$ and $\bz= 1-\frac 1 p>0$.
In this case, $2^{1/p-1}\le B^*/B_*\le 2^{1-1/q}$.
For fixed $p=5/4$, when $q$ varies over $[2,\, 4.25]$, the curves of $H^o$, $B_*$, $B^*$, and $2^{3/2-1/p} B^*$ are given in Figure 2.
The ratio of the upper and lower bounds is increasing in $q$ but no more than $2$.
\begin{center}{\includegraphics[width=11.0cm,height=7.5cm]{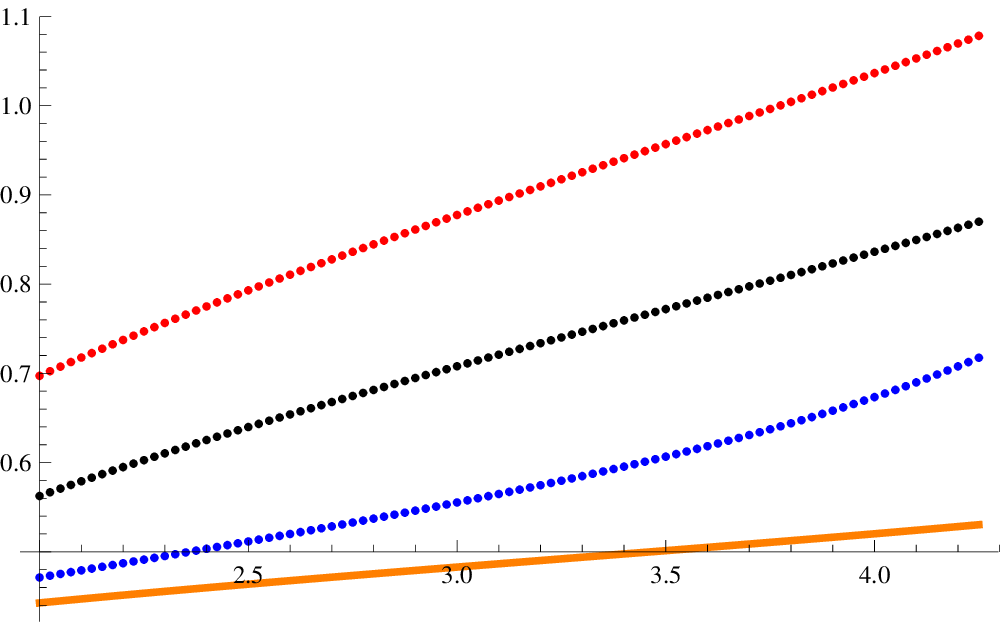}\newline
{\bf Figure 2}\qd The curves from bottom to top are\newline
$H^o$, $B_*$, $B^*$, and $k_{2, p} B^*$, respectively.}\end{center}}
\dexmp

\prf Note that $h\equiv 1$. We have
$$\mu(1, x)\!=\!\int_1^x \!\frac{1}{z^2}=\frac{x-1}{x},\qd
\mu(y, \infty)\!=\!\int_y^\infty\! \frac{1}{z^2}=\frac{1}{y},\qd y(x)\!=\!\frac{x}{x-1},
\qd m(\mu)\!=\!2.$$
Then
$$H_{\mu, \nu}(x, y)=\Big[
\bigg(\frac{x-1}{x}\bigg)^{1/(1-q)}+y^{1/(q-1)} \Big]^{1/q-1}
(y-x)^{(p-1)/p}.$$
Thus,
$$H^{\partial}=\varlimsup_{y\to\infty}H_{\mu, \nu}(x, y)=
\begin{cases}
1 &\text{if $\frac 1 p + \frac 1 q=1$}\\
\infty &\text{if $\frac 1 p + \frac 1 q<1$}\\
0 &\text{if $\frac 1 p + \frac 1 q>1$}.
\end{cases}$$
Hence the inequality does not hold if $\frac 1 p + \frac 1 q<1$.
When $p=2$, this means that the inequality does not hold whenever
$q>2$. The assertion is known as a sharp result for Nash inequality
(cf. \rf{cmf05}{Table 8.2}). Next, we have
$$\aligned
H^o&=2^{1/q-1}\sup_{x\in (1,\,2]}\bigg(\frac{x-1}{x}\bigg)^{1/q} \bigg(\frac{x}{x-1}-x\bigg)^{(p-1)/p}\\
&=2^{1/q-1}\sup_{x\in (1,\,2]}\bigg(1+\frac{1}{x-1}\bigg)^{1-\frac 1 p-\frac 1 q}(2-x)^{1-\frac 1 p}.\endaligned$$
When $\frac 1 p + \frac 1 q=1$, we have $H^o=2^{1/q-1}$, and so $B^*=B_*=H^{\partial}=1$.
When $\frac 1 p + \frac 1 q>1$, we have
$$\aligned
B^*&=\sup_{x<y}\bigg[\bigg(1-\frac 1 x\bigg)^{\frac{p}{(1-p)q}}+y^{\frac{p}{(p-1)q}}\bigg]^{\frac 1 p -1}(y-x)^{1-\frac 1 p},\\
B_*&= \sup_{x<y}\bigg[\bigg(1-\frac 1 x\bigg)^{\frac{1}{1-q}}+y^{\frac{1}{q-1}}\bigg]^{\frac 1 q -1}(y-x)^{1-\frac 1 p}.
\endaligned$$
${H}^o$ has an explicit expression as shown above.
\deprf

\section{Nash inequality, Sobolev-type inequality, and\\ logarithmic Sobolev inequality}

In this section, we study first the Nash inequality and its closely related Sobolev-type
inequality, as a typical application of Theorem \ref{t2-6}. Then we study the logarithmic
Sobolev inequality by the similar method introduced in the paper.

\subsection*{Nash inequality and Sobolev-type inequality}

Recall that the probability measure $\pi$ is defined by $\mu/\mu[-M, N]$. Consider the Nash inequality:
\be \|f-\pi(f)\|_{\mu,\, 2}^{2+4/\gz}\le
A_N \big\|f'\big\|_{\nu,\, 2}^2\|f \|_{\mu,\, 1}^{4/\gz},\qqd \gz>0. \lb{21}\de
It seems more symmetric to replace $\|f\|_{\mu, 1}$ by $\|f -\pi(f)\|_{\mu,\, 1}$ on the right-hand side of (\ref{21}). Let us denote the latter one by (\ref{21})$'$. If (\ref{21}) holds
for every absolutely continuous $f\in L^2(\mu)$, then so does (\ref{21})$'$ regarding $f-\pi(f)$
as a new $f$. Conversely, assume (\ref{21})$'$. Since $\|f -\pi(f)\|_{\mu,\, 1}\le 2 \|f\|_{\mu,\, 1}$, we certainly have (\ref{21}). Hence (\ref{21}) and (\ref{21})$'$ are equivalent.

It is known (cf. \ct{myh02b}, \rf{cmf05}{\S 4.8 and \S 6.5} for related results and more references) that this inequality, when $\gz>2$,  is equivalent to the Sobolev-type inequality
\be \|f-\pi(f)\|_{\mu,\, 2\gz/(\gz -2)}^2\le A_S \big\|f'\big\|_{\nu,\, 2}^2. \lb{22}\de
Now, as an application of Theorem \ref{t2-6}, we have the following result.

\thm\lb{t3-1}{\cms Let $\mu[-M, N]<\infty$ and $\mu_{\text{\rm pp}}=0$. Then
\begin{itemize}
\item[(1)] when $\gz>2$, the Nash inequality (\ref{21}) or equivalently, the Sobolev-type inequality (\ref{22}) holds iff
$H^o\vee H^{\partial}<\infty,$ where $H^o$ and $H^{\partial}$
are defined by (\ref{15}) and (\ref{16}), respectively, with
$p=2$ and $q=2\gz/(\gz-2)$. Furthermore, we have
$$B_{S*}\le A_S\le 4 B_S^*,$$
where
$$\aligned
B_S^{*}&=\sup_{x<y}\Big[\mu(-M, x)^{\frac{2}{\gz}-1}+\mu(y, N)^{\frac{2}{\gz}-1}\Big]^{-1}\,{\hat\nu}[x, y],\\
B_{S*}&=  \inf_{x<y}\Big[\mu(-M, x)^{\frac{4}{2+\gz}-1}+\mu(y, N)^{\frac{4}{2+\gz}-1}\Big]^{-\frac{2}{\gz}-1}\,{\hat\nu}[x, y],
\endaligned
$$
and $B_{S*}\le B_S^{*}\le 2^{2/\gz} B_{S*}$.
\item[(2)] Let $M, N=\infty$. If
$\hat\nu(-\infty, \uz]\wedge \hat\nu [\uz, \infty)=\infty$,
$$\aligned
&\varlimsup_{x\to-\infty} \mu(-\infty, x)\,{\hat\nu}(x, \uz)<\infty \qd
\text{\cms and}\qd \varlimsup_{y\to \infty} \mu(y, \infty)\,{\hat\nu}(\uz, y)<\infty,\\
&\varliminf_{x\to-\infty} \mu(-\infty, x)\,{\hat\nu}(x, \uz)>0\qd
\text{\cms and}\qd \varliminf_{y\to \infty} \mu(y, \infty)\,{\hat\nu}(\uz, y)>0
\endaligned$$
for some $\uz\in {\mathbb R}$, then the Nash inequality (\ref{21}) does not hold whenever $\gz\in (0, 2]$.
\end{itemize}
}
\dethm

\prf (a) The assertion in Part (1) on the estimates of $A_S$
is a straightforward consequence of Theorem \ref{t2-6} with
$p=2$, $q=2\gz/(\gz-2)$, and
$\big(A_S, B_S^*, B_{S*}\big)=\big(A^2, B^{*\,2}, B_*^2\big)$.
The criterion is a copy of Corollary \ref{t2-6-1}.

(b) The restriction ``$\gz>2$'' comes from the reduction of Nash
inequality to the Sobolev-type one. This costs a smaller gap
of the criterion for the Nash inequality,
  marked as $(\vz)$ in the last line on page 15 and the last sentence in Theorem 1.10 of \ct{cmf05}, for instance. The restriction is recently removed in \ct{wj12} in the discrete
situation. Here we show that a direct proof is also possible
under a technical condition.

To see that the Nash inequality (\ref{21}) does not hold for $\gz\in (1, 2]$, rewrite for a moment the inequality as
\begin{equation}\text{\rm Var}_{\mu}(f)^{r}\le A_r \big\|f'\big\|_{\nu,\, 2}^2,\qqd f\in L^2(\mu),\; \|f\|_{\mu,\,1}=1,\label{23}\end{equation}
where $A_r$ denotes the optimal constant.
When $\gz$ varies 0 to 2, $r$ moves from $\infty$ to 2.
By the splitting technique (replacing $\|f\|_{\mu, 1}$ by
$\|f-\pi(f)\|_{\mu, 1}$ on the right-hand side of (\ref{21})), we may consider the half-space $(-M, N)=(0, \infty)$ only, and reduce (\ref{23}) to
$$\|f\|_{\mu,\,2}^{2 r}\le C_r \big\|f'\big\|_{\nu,\, 2}^2,\qqd f(0)=0, \; f\in L^2(\mu),\; \|f\|_{\mu,\,1}=1,   $$
Since $\|f\|_{\mu,\,2}\ge \|f\|_{\mu,\,1}=1$, it is clear that
the last inequality becomes stronger when $r$ increases. Thus, it is sufficient to show
that the inequality (\ref{23}) does not hold when $r=2$ (i.e., $\gz=2$).

To do so, fix  a point $y>0$ and let
$$f(x)={\hat\nu}[0, x\wedge y],\qqd x\ge 0.
$$
Then
$$\aligned
\|f\|_{\mu,\,1}&=\int_0^\infty \mu(\d z)\,{\hat\nu}[0, z\wedge y]
 =\int_0^y h \,\mu(\cdot, \infty),\\
\|f\|_{\mu,\,2}^{2+4/\gz}&\ge \mu(y, \infty)^{1+2/\gz}\,{\hat\nu}[0, y]^{2+4/\gz},\\
\big\|f'\big\|_{\nu,\,2}^2&={\hat\nu}[0, y].
\endaligned
$$
Now, we have
$$\aligned
\frac{\|f\|_{\mu,\,2}^{2+4/\gz}}{\big\|f'\big\|_{\nu,\,2}^2\|f\|_{\mu,\,1}^{4/\gz}}
&\ge \frac{\big(\mu(y, \infty)\,{\hat\nu}[0, y]^2\big)^{1+2/\gz}}{{\hat\nu}[0, y]\big[\int_0^y h\, \mu(\cdot, \infty)\big]^{4/\gz}}\\
&= \frac{\mu(y, \infty)^{1+2/\gz}\,{\hat\nu}[0, y]^{1+4/\gz}}{\big[\int_0^y h\, \mu(\cdot, \infty)\big]^{4/\gz}}\\
&= \frac{\mu(y, \infty)^{1-2/\gz}\,{\hat\nu}[0, y]}{\big[\int_0^y h\, \mu(\cdot, \infty)\big/\big(\mu(y, \infty)\,{\hat\nu}[0, y]\big)+1\big]^{4/\gz}}.
\endaligned$$
Since
$\int_0^y h\, \mu(\cdot, \infty)> \mu(y, \infty)\,{\hat\nu}[0, y]$, when $\gz=2$, we need only to study the ratio
$${\hat\nu}[0, y]\bigg[\frac{\mu(y, \infty)\,{\hat\nu}[0, y]}
{\int_0^y h\, \mu(\cdot, \infty)}\bigg]^{2}.$$
By assumption, if $\int_0^\infty h\, \mu(\cdot, \infty)<\infty$, then the right-hand side goes to infinity as so does $y$ by assumption again. This implies that $A_N=\infty$. Therefore the Nash inequality (\ref{21}) does not hold at $\gz\in (0, 2]$ in this case.

Next, if $\int_0^\infty h\, \mu(\cdot, \infty)=\infty$, then
$$\aligned
\frac{{\hat\nu}[0, y]}{\big[\int_0^y h\, \mu(\cdot, \infty)\big]^{2}}
&\sim \frac{h(y)}{h(y)\mu(y, \infty)\int_0^y h\, \mu(\cdot, \infty)}\qd\text{(by l'H\^ opital's rule)}\\
&=\frac{{\hat\nu}[0, y]}{\mu(y, \infty){\hat\nu}[0, y]\int_0^y h\, \mu(\cdot, \infty)}\\
&\sim\frac{{\hat\nu}[0, y]}{\int_0^y h\, \mu(\cdot, \infty)}\\
&\sim \frac{h(y)}{h(y)\mu(y, \infty)}\\
&=\frac{1}{\mu(y, \infty)}\to \infty \qd\text{as}\qd y\to\infty.
\endaligned$$
Therefore we have arrived the required assertion again.
\deprf

We remark that the condition
$$\varlimsup_{x\to-\infty} \mu(-\infty, x)\,{\hat\nu}(x, \uz)<\infty \qd
\text{and}\qd \varlimsup_{y\to \infty} \mu(y, \infty)\,{\hat\nu}(\uz, y)<\infty$$
in Theorem \ref{t3-1}\,(2) means that the corresponding diffusion
process is exponentially ergodic, otherwise, the Nash inequality
can not hold since the latter inequality is stronger than the former ergodicity (cf. \rf{cmf05}{Table 5.1 and Theorem 1.9}).
Actually, by the cited results, once the transition probability of the process has a density, one can even
assume a stronger condition that $\int_0^\infty h\, \mu(\cdot, \infty)<\infty$, then the proof above can be simplified.

The rough factor $4^{1+1/\gz}\,(\gz>2)$ in Theorem \ref{t3-1}\,(1) is clearly smaller than the first factor $8$ given by
\rf{cmf05}{Theorem 6.8}. The second factor given by the cited theorem is clearly less sharp
than the first one and so is than what we have here.

\subsection*{Logarithmic Sobolev inequality}

We now turn to study the logarithmic Sobolev inequality
\be \ent \big(f^2\big)\le A_{LS} \big\|f'\big\|_{\nu,\, 2}^2,  \lb{24}\de
where
$$\ent (f)=\int_{-M}^N f \log\bigg(\frac{f}{\pi(f)}\bigg)\d\pi
\qd\text{\;for $f\ge 0$}.$$

\thm\lb{t3-2}{\cms Let $\mu[-M, N]<\infty$ and $\mu_{\text{\rm pp}}=0$. Then the optimal constant $A_{LS}$ in (\ref{24})
satisfies the following estimates
$$B_*\le  A_{LS}\le 4 B^*,$$
where
$$\aligned
{B^*}^{-1}&=\inf_{x<y}\bigg\{{\hat\nu}[x, y]^{-1}
\bigg(\bigg[\pi[-M, x]\log\bigg(1+\frac{e^2}{\pi[-M, x]}\bigg)\bigg]^{-1}\\
&\qqd\qqd\qqd\qqd\qd +\bigg[\pi[y, N]\log\bigg(1+\frac{e^2}{\pi[y, N]}\bigg)\bigg]^{-1}\bigg)\bigg\},\\
B_*^{-1}&=\inf_{\uz\in (-M, N)}\inf_{(x, y)\ni \uz}
\bigg\{{\hat\nu}[x, y]^{-1}\bigg(\bigg[\pi[-M, x]\log\bigg(1+\frac{1-\pi[-M,\, \uz]}{\pi[-M, x]}\bigg)\bigg]^{-1}\\
&\hspace{13em}+\bigg[\pi[y, N]\log\bigg(1+\frac{1-\pi[\uz, N]}{\pi[y, N]}\bigg)
\bigg]^{-1}\bigg)\bigg\},
\endaligned$$
or alternatively,
\begin{align}
B_*^{-1}&=\inf_{x< y}\bigg\{{\hat\nu}[x, y]^{-1}\bigg(\bigg[\pi[-M, x]\log\bigg(1+\frac{z^*(x, y)}{\pi[-M, x]}\bigg)\bigg]^{-1}\nnb\\
&\hspace{7em}\;+\bigg[\pi[y, N]\log\bigg(1+\frac{1-z^*(x, y)}{\pi[y, N]}\bigg)\bigg]^{-1}
\bigg),\lb{28-1}
\end{align}
where $z^*(x, y)$ is the unique solution to the equation
\begin{align}
&\bigg[\pi[-M, x]\log\bigg(1+\frac{z}{\pi[-M, x]}\bigg)\bigg]^2
\bigg(1+\frac{z}{\pi[-M, x]}\bigg)\nnb\\
&\qd =\bigg[\pi[y, N]\log\bigg(1+\frac{1-z}{\pi[y, N]}\bigg)\bigg]^2
\bigg(1+\frac{1-z}{\pi[y, N]}\bigg),\qqd z\in (0, 1).\lb{28-2}
\end{align}
In particular, we have
$$\aligned
B_*^{-1}&\le
\inf_{(x,\, y)\ni m(\pi)}\bigg\{{\hat\nu}[x, y]^{-1}
\bigg(\bigg[\pi[-M, x]\log\bigg(1+\frac{1}{2\pi[-M, x]}\bigg)\bigg]^{-1}\\
&\qqd\qqd\qqd\qqd\qqd\qd +\bigg[\pi[y, N]\log\bigg(1+\frac{1}{2\pi[y, N]}\bigg)\bigg]^{-1}\bigg)\bigg\},
\endaligned$$
where $m(\pi)$ is the median of $\pi$.}
\dethm

\prf
(a) Upper bound. Even though this theorem is not a consequence of Theorem \ref{t2-9}, the
idea of the proof for the upper bound is more or less the same as we used several times
in the last two sections. Given $f\in L^2(\pi)$ with $\nu\big({f'}^2\big)\in (0, \infty)$ and $\uz\in (-M, N)$, let
$$\tilde f=f-f(\uz), \qqd {\tilde f}_-={\tilde f}\mathbbold{1}_{[-M,\, \uz]},\qqd
{\tilde f}_+={\tilde f}\mathbbold{1}_{(\uz,\, N]}.$$
The following facts were proved in the first part of \rf{bfrc}{Proof of Theorem 3} for the specific $\uz=m(\pi)$, but the proof remains true for general $\uz$:
$$\aligned
\ent \big(f^2\big)&\le \ent \big(\tilde f^2\big)+ 2 \pi\big(\tilde f^2\big)\qd\text{(by \rf{cmf05}{Lemma 4.14})}\\
&\le \ent \big(\tilde f_-^2\big)+ 2 \pi\big(\tilde f_-^2\big)+\ent \big(\tilde f_+^2\big)+ 2 \pi\big(\tilde f_+^2\big)
\endaligned
$$
(since $\ent$ is sub-additive: $\ent(f+g)\le \ent(f)+\ent(g) $) and
$$\ent \big(\tilde f_{\pm}^2\big)+ 2 \pi\big(\tilde f_{\pm}^2\big)\le 4 B_{\uz}^{\pm} \nu\big({{\tilde f}_{\pm}}^{\prime\,2}\big),$$
where
$$\aligned
B_{\uz}^-&=\sup_{x<\uz}\pi[-M, x]\log\bigg(1+\frac{e^2}{\pi[-M, x]}\bigg)\,{\hat\nu}[x, \uz],\\
B_{\uz}^+&=\sup_{y>\uz}\pi[y, N]\log\bigg(1+\frac{e^2}{\pi[y, N]}\bigg)\,{\hat\nu}[\uz, y].
\endaligned$$
Thus, we have
$$\aligned
\frac{\ent\big(f^2\big)}{\nu \big({f'}^2\big)}
&\le \frac{\ent\big(\tilde f^2\big)+2 \pi\big(\tilde f^2\big)}{\nu \big({f'}^2\big)}\\
&\le \frac{\ent\big(\tilde f_-^2\big)+2 \pi\big(\tilde f_-^2\big)+\ent\big(\tilde f_+^2\big)+2 \pi\big(\tilde f_+^2\big)}{\nu \big({f_-'}^2\big)+\nu \big({f_+'}^2\big)}\\
&\le \frac{\ent\big(\tilde f_-^2\big)+2 \pi\big(\tilde f_-^2\big)}{\nu \big({f_-'}^2\big)}\bigvee
\frac{\ent\big(\tilde f_+^2\big)+2 \pi\big(\tilde f_+^2\big)}{\nu \big({f_+'}^2\big)}\\
&\le \big(4 B_{\uz}^-\big)\vee \big(4 B_{\uz}^+\big).
\endaligned$$
The original proof for the upper estimate stopped here with $\uz=m(\pi)$.
Because $\uz$ is arbitrary, we obtain
$$\frac{\ent\big(f^2\big)}{\nu \big({f'}^2\big)}\le 4 \inf_{\uz}\big(B_{\uz}^- \vee  B_{\uz}^+\big).$$
Note that the right-hand side is independent of $f$.
By choosing $\bar\uz$ such that $B_{\bar\uz}^- =  B_{\bar\uz}^+$, it follows that
$$\frac{\ent\big(f^2\big)}{\nu \big({f'}^2\big)}\le 4 B_{\bar\uz}^- .$$
Now, since $f$ with $\nu\big({f'}^2\big)\in (0, \infty)$ is arbitrary, we obtain
$$A_{LS}=\sup_{\nu({f'}^2)\in (0, \infty)}\frac{\ent\big(f^2\big)}{\nu \big({f'}^2\big)}\le 4 B_{\bar\uz}^- .$$
Next, as an application of Lemma \ref{t1-2}, we have
$$\aligned
B^*&=\sup_{x<y}\frac{{\hat\nu}\big[x, \bar\uz\,\big] + {\hat\nu}\big(\bar\uz, y\big]}{\big[\pi[-M, x]\log\big(1+\frac{e^2}{\pi[-M, x]}\big)\big]^{-1}+
\big[\pi[y, N]\log\big(1+\frac{e^2}{\pi[y, N]}\big)\big]^{-1}}\\
&\ge B_{\bar\uz}^-\wedge B_{\bar\uz}^+\\
&= B_{\bar\uz}^-.
\endaligned$$
Combining the last two estimates together, we obtain the required upper bound. Once again, the unknown $\bar\uz$ disappears in the
expression of $B^*$.

(b) Lower bound. We adopt a similar method as used in the proof of Lemma \ref{t2-4}. Define
$$\aligned
f(z)=- \mathbbold{1}_{\{z\le \uz\}}\,{\hat\nu}[x\vee z,\, \uz]
+ \mathbbold{1}_{\{z> \uz\}}\,{\hat\nu} [\uz,\, y\wedge z],\qqd z\in [-M, N],\endaligned$$
where $x, y, \uz$ with $(x, y)\ni\uz$ are fixed. First, let us apply \rf{bfrc}{Proof of Theorem 3} to this specific test
function $f$.
$$\aligned
\ent\big(f^2 \mathbbold{1}_{\{z\le \uz\}}\big)
&=\sup\bigg\{\int_{-M}^{\uz} f^2 g\d\pi: \int_{-M}^{\uz} e^g\d\pi\le 1\bigg\}\\
&\ge \sup\bigg\{\int_{-M}^{\uz} f^2 g\d\pi: g\ge 0\; \text{ and }\int_{-M}^{\uz} e^g\d\pi\le 1\bigg\}\\
&\ge {\hat\nu}[x, \uz]^2
\sup\bigg\{\!\int_{-M}^{\uz}\! \mathbbold{1}_{[-M,\, x]}\, g\d\pi: g\ge 0 \;\text{ and }\!\int_{-M}^{\uz}\!\! e^g\d\pi\le 1\bigg\}.
\endaligned$$
Applying \rf{bfrc}{Lemma 6} to the last supremum, it follows that
\begin{align}
\ent\big(f^2 \mathbbold{1}_{[-M,\, \uz]}\big)
&\ge {\hat\nu}[x, \uz]^2 \,\fz(x, \uz),\nnb\\
\fz(x, \uz):&=\pi[-M, x]\log\bigg(1+\frac{1-\pi[-M,\, \uz]}{\pi[-M, x]}\bigg).\lb{26}\end{align}
Symmetrically, we have
\begin{align}
\ent\big(f^2 \mathbbold{1}_{[\uz, N]}\big)
&\ge {\hat\nu}[\uz, y]^2 \,\qz (\uz, y),\nnb\\
\qz (\uz, y):&=\pi[y, N]\log\bigg(1+\frac{1-\pi[\uz, N]}{\pi[y, N]}\bigg).\lb{27}\end{align}

Next, by logarithmic Sobolev inequality,
$$A_{LS}\,\nu\big({f'}^2  \mathbbold{1}_{[-M,\, \uz]}\big)
\ge \ent\big(f^2 \mathbbold{1}_{[-M,\, \uz]}\big),$$
it follows that
$$A_{LS}\,\nu\big({f'}^2  \mathbbold{1}_{[-M,\, \uz]}\big)
\ge {\hat\nu}[x, \uz]^2 \fz (x, \uz).$$
Similarly,
$$A_{LS}\,\nu\big({f'}^2  \mathbbold{1}_{[\uz, N]}\big)
\ge {\hat\nu}[\uz, y]^2 \,\qz(\uz, y).$$
We now arrive at the place different from the known proofs. Summing up the last two inequalities, it follows that
$$A_{LS}\,\nu\big({f'}^2\big)
\ge {\hat\nu}[x, \uz]^2 \,\fz(x, \uz)+{\hat\nu}[\uz, y]^2 \,\qz(\uz, y).$$
Since one can replace $\nu\big({f'}^2\big)$ by ${\hat\nu}\big({f'}^2\big)$ in the original inequality, by definition of $f$, we have ${\hat\nu}\big({f'}^2\big)={\hat\nu}[x, y]$, and so
$$A_{LS}\ge \frac{1}{{\hat\nu}[x, y]}
\Big[{\hat\nu}[x, \uz]^2 \fz(x, \uz)+{\hat\nu}[\uz, y]^2 \,\qz(\uz, y)\Big].$$
Therefore, we have
$$A_{LS}\ge
{\hat\nu}[x, y]\bigg[\bigg(\frac{{\hat\nu}[x, \uz]}{{\hat\nu}[x, y]}\bigg)^2\fz(x, \uz)+
\bigg(1-\frac{{\hat\nu}[x, \uz]}{{\hat\nu}[x, y]}\bigg)^2\qz(\uz, y)\bigg].$$
Noting that the function $c_1 z^2+ c_2 (1-z)^2$ on $[0, 1]$ achieves
its minimum $\big(c_1^{-1}+c_2^{-1}\big)^{-1}$ at $z^*=c_2/(c_1+c_2)$, it follows that
$$A_{LS}\ge
{\hat\nu}[x, y]\big(\fz(x, \uz)^{-1}+
\qz(\uz, y)^{-1}\big)^{-1}.$$
Since $(x, y)\ni \uz$ are arbitrary, we finally arrive at
\begin{align}
A_{LS}\ge\sup_{\uz\in (-M, N)}\sup_{(x, y)\ni \uz}
{\hat\nu}[x, y]\big(\fz(x, \uz)^{-1}+
\qz(\uz, y)^{-1}\big)^{-1}.\end{align}
This gives us the first version of $B_*$.
Then the final assertion of the theorem follows by setting $\uz=m(\pi)$.

(c) We now prove the alternative assertion of $B_*$. Since supremums are exchangeable, one may rewrite $B_*$ as
$$B_*=\sup_{x< y}\bigg\{{\hat\nu}[x, y]\sup_{\uz\in (x, y)}
\big(\fz(x, \uz)^{-1}+
\qz(\uz, y)^{-1}\big)^{-1}\bigg\}.$$
Fix $x<y$ and make a change of the variable $\uz$ by $z=1-\pi[-M, \uz]$. Then $1-\pi[\uz, N]=1-z$ and the functions $\fz(x, \uz)$ and $\qz(\uz, y)$ become ${\tilde\fz}(x, z)$ and ${\tilde\qz}(z, y)$, respectively.
The supremum above should be achieved at the point for which the derivative in $z$ of
$\big({\tilde\fz}(x, z)^{-1}+ {\tilde\qz}(z, y)^{-1}\big)^{-1}$ vanishes.  This leads to the unique solution $z^*=z^*(x, y)$ to equation
(\ref{28-2}). Then we obtain (\ref{28-1}).
\deprf

We mention that there is a large number of publications on the logarithmic Sobolev inequalities,
in the one-dimensional case for instance,
one may refer to \ct{myh02a}, \ct{bfrc}, \rf{cmf05}{\S 4.6 and \S 6.6} for related results and more references. Generally speaking,
Theorem \ref{t3-2} clearly improves \rf{bfrc}{Theorem 3}
(since $\az\vee\bz$ is replaced by $\az+\bz$), to which
the factor of the upper and lower bounds is at most 16, the best one we have known up to now. In the special case that the measures
$\pi$ and $\hat\nu$ are symmetric with respect to $m(\pi)$, the computation of $B^*$ and $B_*$ can be reduced to half space. Then
Theorem \ref{t3-2} coincides with \rf{bfrc}{Theorem 3}. Besides,
having Theorem \ref{t3-2} at hand, it should be easy to introduce
the corresponding $H^o$ and $H^{\partial}$ as we did in the previous sections.

The methods introduced in the paper should have more applications. For instance, in parallel to the proof of Theorem \ref{t3-2}, we may have an improved version of the Sobolev inequality and the
Lata{\l}a--Oleszkiewicz inequality presented by \rf{bfrc}{Theorems 11 and 13}.
\bigskip

\nnd{\bf Acknowledgments}. {\small
The results of the paper have been presented in our seminar.
Thanks are given to Y.H. Mao, F.Y. Wang, Y.H. Zhang, and the participants for their
helpful comments and suggestions which lead to some improvements
of the paper. Thanks are also given to a referee for correcting
a number of typos in a previous version of the paper.
Research supported in part by the
         National Natural Science Foundation of China (No. 11131003),
         and by the ``985'' project from the Ministry of Education in China.
}

\nnd {\small School of Mathematical Sciences, Beijing Normal University,\newline
Laboratory of Mathematics and Complex Systems (Beijing Normal University),\newline
\mbox{\qqd} Ministry of Education, Beijing 100875,
    The People's Republic of China.\newline E-mail: mfchen@bnu.edu.cn\newline Home page:
    http://math.bnu.edu.cn/\~{}chenmf/main$\_$eng.htm
}

\end{document}